\documentclass{article}

\usepackage{amssymb,amsmath,url,pgf,graphicx,subfig}
\usepackage[ruled,linesnumbered,longend]{algorithm2e}

\newcommand{\pf}{Proof. }
\newcommand{\diag}{\hbox{diag}}
\newcommand{\qed}{\hfill $\square$\smallskip}
\newcommand{\fl}{\hbox{fl}}

\newcommand{\C}{\mathbb{C}}
\newcommand{\CC}{\mathcal{C}}
\newcommand{\TT}{\mathcal{T}}
\newcommand{\LL}{\mathcal{L}}
\newcommand{\cu}{{\tt i}}
\newcommand{\wh}{\widehat}
\newcommand{\wt}{\widetilde}
\newcommand{\D}{\mathcal D}

\newcommand{\uno}{{\bf 1}}
\newcommand{\dotleq}{\begin{array}{c}^{\displaystyle \cdot}\\[-2ex]\leq\end{array}}

\newtheorem{theorem}{Theorem}
\newtheorem{proposition}[theorem]{Proposition}
\newtheorem{remark}[theorem]{Remark}

\newcommand{\ud}{\mathrm{d}}
\newcommand{\scal}{\mathcal{S}}

\usepackage{color}

\definecolor{darkmagenta}{rgb}{0.5,0,0.5}
\definecolor{darkgreen}{rgb}{0,0.6,0}
\definecolor{darkblue}{rgb}{0,0,0.6}
\definecolor{darkred}{rgb}{0.8,0,0}
\definecolor{mellow}{rgb}{.847, 0.72, 0.525}

\title{Computing the Exponential of Large Block-Triangular Block-Toeplitz
  Matrices Encountered in Fluid Queues}

\author{D.A. Bini, Universit\`a di Pisa\\
S. Dendievel, Universit\'e libre de Bruxelles\\
G. Latouche, Universit\'e libre de Bruxelles\\
B. Meini,  Universit\`a di Pisa}

\begin{document}
\maketitle

\begin{abstract}
The Erlangian approximation of Markovian fluid queues leads to the
problem of computing the matrix exponential of a subgenerator having a
block-triangular, block-Toeplitz structure. 
To this end, we propose some algorithms which exploit the Toeplitz
structure and the properties of generators. Such algorithms allow to
compute the exponential of very large matrices, which would otherwise
be untreatable with standard methods. We also prove interesting decay properties of the exponential of a generator having a
block-triangular, block-Toeplitz structure.
\end{abstract}


\section{Introduction}
The problem we consider here is to compute the exponential of an upper block-triangular,
block-Toeplitz matrix, that is, a matrix of the
kind
\begin{equation}\label{U}
\mathcal T(U)=\left[\begin{array}{cccc}
U_0 & U_1 & \ldots & U_{n-1}\\
    & U_0 & \ddots & \vdots \\
    &     & \ddots & U_1 \\
0    &     &        & U_0
\end{array}\right],
\end{equation}
where $U_i$, $i=0,\ldots,n-1$, are $m\times m$ matrices.
Our interest stems from the analysis in Dendievel and
Latouche~\cite{dl14} of the Erlangization method for
Markovian fluid models, but the story goes further back in time.

\subsection{Origin of the problem}

The Erlangian approximation method was introduced in Asmussen {\it et
  al.} \cite{aau02} in the context of risk processes; it was picked up
in Stanford {\it et al.} \cite{sabbdl05} where a connection is
established with fluid queues.  Other relevant references
are Stanford {\it et al.} \cite{slwbh05} where the focus is on
modelling the spread of forest fires, and Ramaswami {\it et al.}
\cite{rws08} where some basic algorithms are developed.

Markovian fluid models are two-dimensional processes $\{(X(t),
\varphi(t)): t \in \mathbb{R}^+\}$ where $\{\varphi(t)\}$ is a Markov
process with infinitesimal generator $A$ on the state space $\{1,
\ldots , m\}$; to each state $i$ is associated a rate of growth $c_i
\in \mathbb{R}$ and $X(t)$ is controlled by $\varphi(t)$ through
the equation 
\[
X(t) = X(0) + \int_0^t c_{\varphi(s)} \, \ud s, \qquad \mbox{for $t
  \geq 0$.}
\]
Performance measures of interest include the distributions of $X(t)$
and of various first passage times.  
 Usually, $\varphi(t)$ is called the phase of the process at time $t$
and $X(t)$ its level, and  the phase space $\{1, \ldots,
m\}$ is partitioned into three subsets $\scal_+$, $\scal_-$ and $\scal_0$ such that $c_i >0$,
$c_i < 0$ or $c_i =0$ if $i$ is in $\scal_+$, $\scal_-$ or $\scal_0$,
respectively.   To simplify our presentation without missing any
important feature, we assume below that $\scal_0$ is empty.

The first return probabilities of
$X(t)$ to its initial level $X(0)$ play a central role in the analysis
of fluid queues.  It is
customary to define two matrices $\Psi$ and $\hat \Psi$ of first
return probabilities:
\[
\Psi_{ij} = \Pr[\tau < \infty, \varphi(\tau) = j | X(0)=0,
\varphi(0)=i],   \qquad \mbox{$i \in \scal_+$, $j \in \scal_-$,}
\]
and
\[
\hat{\Psi}_{ij} = \Pr[\tau < \infty, \varphi(\tau) = j | X(0)=0,
\varphi(0)=i],   \qquad \mbox{$i \in \scal_-$, $j \in \scal_+$,}
\]
where $\tau = \inf\{t > 0: X(t) =0\}$ is the first passage time to
level 0.  Thus, the entries of $\Psi$ and $\hat\Psi$ are the
probability of returning  to the initial level after having started in
the upward, and the downward directions, respectively.

If the process starts from some level $x > 0$, then
\[
\Pr[\tau < \infty, \varphi(\tau) = j | X(0)=x, \varphi(0)=i] 
=
(\begin{bmatrix} I \\ \Psi \end{bmatrix} e^{Hx})_{ij}
\qquad \mbox{$i \in \{1, \ldots, m\}$, $j \in \scal_-$;}
\]
here, $H$ is a square matrix on $\scal_- \times \scal_-$ and is given by
\[
H = |C_-|^{-1} A_{--}  +   |C_-|^{-1} A_{-+}  \Psi,
\]
where $A_{--}$ and    $A_{-+}$ are submatrices of the generator $A$,
indexed by  $\scal_- \times \scal_-$ and $\scal_- \times \scal_+$, respectively,
and $|C_-|$ is a diagonal matrix with $|c_i|, i \in \scal_-$  on the
diagonal.
A similar equation holds for $x < 0$.
The matrices $\Psi$ and $\hat \Psi$ are solutions of algebraic Riccati
equations and their resolution has been the object of much attention.
Very efficient algorithms are available, and we refer to Bini {\it et
  al.}~\cite{bim12} and Bean {\it et al.}~\cite{brt05b}.

The Erlangian approximation method is introduced in \cite{aau02} to
determine the detailed distribution of $\tau$.  The idea is that, to
compute the probability
\[
F(t_0; i, x) = \Pr[\tau < t_0 | X(0)=x, \varphi(0)=i],
\qquad \mbox{$x>0$, $i \in \{1, \ldots , m\}$}
\]
for a fixed value $t_0$, it is convenient to replace $t_0$ by a random
variable $T$ with an Erlang distribution, with parameters $(n/t_0, n)$
for some positive integer $n$.  The random variable $T$ has
expectation $t_0$ and variance $t_0/n$, so that $F(T;i,x)$ is a good
approximation of $F(t_0;i,x)$ if $n$ is large enough.  From a
computational point of view, the advantage is that one replaces
systems of integro-differential equations by linear equations.

The long and the short of it is that the original system is replaced by the process
$\{(X(t), \Phi(t))\}$ with a two-dimensional phase $\Phi(t)=(\beta(t),
\varphi(t))$ on the state space $\{1, \ldots, n\} \times \{1, \ldots,
m\} \cup\{0\}$ and with the generator
\[
Q =
\begin{bmatrix}
A-\nu I & \nu I & & 0 & 0 \\
 & A-\nu I & \ddots \\
& & \ddots & \nu I \\
 & & &    A-\nu I & \nu \uno \\
0 & & & 0 & 0
\end{bmatrix},
\]
where $\nu = n/t_0$.  The physical interpretation is that the
absorbing state 0 is entered at the random time $T$, and the component
$\beta$ of the phase marks the progress of time towards $T$.  Some authors (for
instance \cite{aau02, sabbdl05}) report that good approximations may be
obtained with small values of $n$.

Because of the Toeplitz-like structure of $Q$, the matrices $\Psi$ and
$H$ are both upper block-triangular block-Toeplitz and it is
interesting to use the Toeplitz structure in order to reduce the cost
when $n$ is large.  This is done in \cite{rws08} for the matrix
$\Psi$.  Here we address the question of efficiently computing the
exponential matrix $e^{Hx}$ for a given value of $x$, where $H$ has the structure of \eqref{U}.  We shall assume
without loss of generality that $x=1$.

\subsection{Main results}

We recall that the exponential function can be extended to a matrix variable by defining
\begin{equation}\label{eq:mexpf}
e^X=\sum_{i=0}^\infty \frac 1{i!}X^i.
\end{equation}
For more details on the matrix exponential and more generally on matrix functions we refer the reader to Higham~\cite{highambook}.

The matrix  $\TT(U)$ defined in \eqref{U} is of order $nm$ and it may be huge,
since a larger $n$ leads to a better Erlangian approximation, while
the size $m$ of the blocks is generally small.  The matrix $\TT(U)$ is a
subgenerator, i.e., it has negative diagonal entries, nonnegative
off-diagonal entries, and the sum of the entries on each row is
nonpositive. 

Since block-triangular block-Toeplitz matrices are closed under matrix multiplication,  it follows from \eqref{eq:mexpf} that 
the matrix exponential $e^{\TT(U)}$ is also an upper block
triangular, block-Toeplitz matrix; in particular, the diagonal blocks
of $e^{\TT(U)}$ coincide with $e^{U_0}$. Moreover, it is known that 
the matrix
$e^{\TT(U)}$ is nonnegative and substochastic.

The problem of the computation of the exponential of a generator has
been considered in Xue and Ye~\cite{xy08,xy13} and by Shao et
al. \cite{shao}, where the authors propose component-wise accurate
algorithms for the computation. These algorithms are efficient for
matrices of small size. For the Erlangian approximation problem, these
algorithms are useless for the large size of the matrices involved.
Recently, some attention has been given to the computation of the
exponential of general Toeplitz matrices by using Arnoldi 
method (Lee {\it et al.}~\cite{LeePangSun}, Pang and
Sun~\cite{PangSun}).  

In our framework, Toeplitz matrices are block-triangular so that they
form a matrix algebra. This property is particularly effective for the
design of efficient algorithms and we propose some numerical methods
that exploit the
block-triangular block-Toeplitz structure and the generator
properties. Unlike the general methods, our algorithms allow one to deal
with matrices $\TT(U)$ of very large size.

Two methods rely on spectral and computational properties of
block-circulant and block $\epsilon$-circulant matrices
(Bini~\cite{bini84}, Bini {\it et al.}~\cite{blmbook}) and on the use
of Fast Fourier Transforms
(FFT). Recall that block $\epsilon$-circulant matrices have the form
\[
\mathcal C_\epsilon(U)=\left[\begin{array}{cccc}
U_0 & U_1 & \ldots & U_{n-1}\\
\epsilon U_{n-1}    & U_0 & \ddots & \vdots \\
\vdots    &  \ddots   & \ddots & U_1 \\
 \epsilon U_{1}   &  \ldots   &   \epsilon U_{n-1}     & U_0
\end{array}\right],
\] 
and that a block-circulant matrix is a block $\epsilon$-circulant
matrix with $\epsilon=1$. For simplicity, we denote by $\mathcal C(U)$ the 
block $1$-circulant matrix $\mathcal C_1(U)$.

Since block $\epsilon$-circulant matrices can be block-diagonalized by
FFT \cite{bini84}, the computation of the exponential of an $n\times
n$ block $\epsilon$-circulant matrix with $m\times m$ blocks can be
reduced to the computation of $n$ exponentials of $m\times m$
matrices. These latter exponentials are independent from each other
and can be computed simultaneously with a multi-core architecture at
the cost of a single exponential.

The idea of the first method is to approximate $e^{\TT(U)}$ by
$e^{\CC_\epsilon(U)}$ where $\epsilon\in\mathbb C$ and $|\epsilon|$ is
sufficiently small.  We analyse the error and are thereby able to
choose the value of $\epsilon$ which gives a good balance between the
roundoff error and the approximation error.  In fact, the
approximation error grows as $O(\epsilon)$ while the roundoff error is
$O(\mu\epsilon^{-1})$, where $\mu$ is the machine precision. This
leads to an overall error which is $O(\mu^{1/2})$. By using the fact
that the solution is real, by choosing $\epsilon$ a pure imaginary
number we get an approximation error $O(\epsilon^2)$ which leads to an
overall error $O(\mu^{2/3})$.

Since the approximation error is a power series in $\epsilon$, we
devise a further technique which consists in averaging the solutions
computed with $k$ different values of $\epsilon$. This way, we are
able to cancel out the components of the error of degree less than
$\epsilon^{2k}$. This leads to a substantial improvement of the
precision. Moreover, since the different computations are independent
from each other, the computational cost in a multicore architecture is
independent of $k$.

In our second approach, the matrix $\TT(U)$ is embedded into a
$K\times K$ block-circulant matrix $\CC(U^{(K)})$, where $K$ is
sufficiently large, and an approximation of ${\rm e}^{\TT(U)}$ is
obtained from a suitable submatrix of $e^{\CC(U^{(K)})}$.  The computation
of $e^{\CC(U^{(K)})}$ is reduced to the computation of $K$ exponentials of
$m\times m$ matrices, and our error analysis allows one to choose the
value of $K$ so as to guarantee a given error bound in the computed
approximation.

The third numerical method consists in specializing the shifting and
Taylor series method of \cite{xy13}. The block-triangular Toeplitz
structure is exploited in the FFT-based matrix multiplications
involved in the algorithm, leading to a reduction of the computational
cost. The algorithm obtained in this case does not seem well suited
for an implementation in a multicore architecture.

We compare the three numerical methods, from a theoretical as well as
from a numerical point of view.  From our analysis, we conclude that
the method based on $\epsilon$-circulant matrices is the fastest and
provides a reasonable approximation to the solution. Moreover, by
applying the averaging technique we can dramatically improve the
accuracy.
The method based on embedding and the one based on power series
perform an accurate computation but are slightly more expensive.  

It must be emphasised that the use of FFT makes the algorithms
norm-wise stable but that component-wise stability is not
guaranteed. In consequence, the matrix elements with values of modulus
below the machine precision may not be well approximated in terms of
relative error.

The paper is organised as follows. In Sections \ref{sec:der} and
\ref{sec:fastcomp}, we recall properties of the exponential of a
subgenerator and of its derivatives, and some basic properties of
block-Toeplitz and block-circulant matrices which are used in our
algorithms. In Section \ref{sec:expcirc}, we show how to compute the
exponential of a block $\epsilon$-circulant matrix by using fast
arithmetic based on FFT and we perform an error analysis.  We present
in Section \ref{sec:exptt} the algorithms to compute the exponential
of $\TT(U)$: first we analyse the decay of off-diagonal entries of the
matrix exponential, next we describe the new methods and perform an
error analysis. We conclude with numerical experiments in Section
\ref{sec:exper}.

\section{The exponential of a subgenerator and its\\  derivatives }\label{sec:der}

\subsection{The exponential of a subgenerator}
A subgenerator of a Markov process is a matrix $Q$ of real numbers such that
the off--diagonal entries of $Q$ are nonnegative, the diagonal entries 
are negative, and the sum of the entries on each row is nonpositive.  
We denote by $\uno$ the column vector with all entries equal to 1,
with size according to the context. If $Q$ is a subgenerator, then
$Q\uno\le 0$ and $Q$ is called a generator if the row sum on all rows
is zero.

Let $\sigma=\max_i(-q_{ii})$. The matrix $V=Q+\sigma
I$ is a nonnegative matrix, and we may write $e^Q=e^{V-\sigma
  I}=e^{-\sigma}e^V$.  From the latter equality it follows that the
matrix exponential $e^Q$ is nonnegative. Moreover, since $Q\uno\le 0$
it follows that $V\uno=Q\uno+\sigma\uno\le \sigma\uno$. Therefore, in view of \eqref{eq:mexpf},
$e^Q\uno=e^{-\sigma}e^V\uno= e^{-\sigma}\sum_{i=0}^{\infty}\frac
1{i!}V^i\uno\le e^{-\sigma}e^\sigma\uno$. Thus we may conclude that
$e^Q\uno\le\uno$, that is, $e^Q$ is a substochastic matrix.

\subsection{Derivatives and perturbation results}
We recall the definition and some properties of the  G\^{a}teaux and
Fr\'echet derivatives, and their expression for the matrix exponential
function, together with some properties when the matrix is a
subgenerator. We refer the reader to \cite{highambook} for
more details.

The Fr\'echet derivative of a matrix function $f:\C^{n\times
  n}\rightarrow \C^{n\times n}$ at a point $X\in\C^{n\times n}$ along
the direction $E\in\C^{n\times n}$ is the linear mapping $L(X,E)$ in
the variable $E$ such that
\begin{equation}\label{fre}
f(X+E)-f(X)-L(X,E)=o(\| E\|).
\end{equation}
The G\^{a}teaux (or directional)
derivative of  $f:\C^{n\times n}\rightarrow \C^{n\times n}$
at a point $X\in\C^{n\times n}$ along the direction $E\in\C^{n\times n}$ is
\begin{equation}\label{gat}
G(X,E)= \lim_{h\to0}\frac{f(X+hE)-f(X)}{h}.
\end{equation}
If the Fr\'echet derivative exists, then it is equal to the G\^{a}teaux derivative
(\cite[Section 3.2]{highambook}). Such is the case for the matrix
exponential function and we may, therefore, use either definition
\eqref{fre} or \eqref{gat}, depending on which is more convenient; 
we will use the G\^{a}teaux derivative.  From
\cite{nh95},
\begin{equation}\label{frex}
G(tX,E)=\int_0^t e^{X(t-s)} E e^{Xs} ds
\end{equation}
and the following equation gives an expression for the matrix
exponential in terms of G\^{a}teaux derivatives:
\begin{equation}\label{tayl}
e^{t(X+hE)}=\sum_{j=0}^\infty \frac{h^j}{j!}G^{[j]}(tX,E) 
\end{equation}
where we denote by $G^{[j]}(tX,E)$  
the $j$-th  G\^{a}teaux derivative of the matrix function
$e^{tX}$ in the direction $E$, obtained by the recurrence equation
\begin{equation}\label{gatint}
G^{[j]}(tX,E)=j\int_0^t e^{(t-s)X} E G^{[j-1]}(sX,E) ds,~~~j=1,2,\ldots,
\end{equation}
and $G^{[0]}(tX,E)=e^{tX}$.

Recall that if $X$ is a subgenerator, then $e^{tX}$ is a
substochastic matrix for any $t\ge 0$ and in particular $\| e^{tX}
\|_\infty\le 1$.  Therefore, by taking norms in \eqref{frex}, we
obtain the upper bound
\begin{equation}\label{frecbound}
\| G(tX,E)\|_\infty \le \int_0^t \|e^{X(t-s)}\|_\infty \|E\|_\infty
\|e^{Xs}\|_\infty ds\le t \|E\|_\infty
\end{equation}
which may be extended to the $j$-th order G\^{a}teaux derivative as in
the next proposition.

\begin{proposition}\label{prop1}
If $X$ is a subgenerator, then
\begin{equation}\label{boundder}
\| G^{[j]}(tX,E)\|_\infty \le t^j \| E\|_\infty^j 
\end{equation}
for $t \geq 0$, for any $j\ge 0$. Moreover, if $E$ is a nonnegative
matrix, then $ G^{[j]}(tX,E)$ is nonnegative for any $j\ge 0$.
\end{proposition}

\pf
Since the matrix $X$ is a subgenerator, the
matrix $e^{\tau X}$ is nonnegative and substochastic for any $\tau\ge 0$,
therefore $\| e^{\tau X}\|_\infty\le 1$ for any $\tau\ge 0$. 
By using this property, the inequality
\eqref{boundder} can be proved by induction. If $j=0$, then  $\|
G^{[0]}(tX,E)\|_\infty \le 1$. The inductive step is immediately proved,
since from \eqref{gatint} we have
\[
\begin{split}
\|G^{[j]}(tX,E)\|_\infty & \le j\int_0^t \|e^{(t-s)X}\|_\infty \|E\|_\infty
\|G^{[j-1]}(sX,E)\|_\infty ds\\
& \le j \int_0^t \|E\|_\infty^j s^{j-1} ds=\|E\|_\infty^j t^j,
\end{split}
\]
where the last inequality follows from the inductive assumption. If the matrix
$E$ is nonnegative, from the recurrence \eqref{gatint} and from the fact that
 $e^{\tau X}$ is nonnegative and substochastic for any $\tau\ge 0$, it follows
 by induction that $G^{[j]}(tX,E)$ is nonnegative for any $j\ge 0$.
\qed

The following result provides some bounds related to the exponential of
the matrix $\TT(U)$ of \eqref{U} and to its G\^{a}teaux derivative; it will be
used in the next sections to analyse the stability of the algorithm in Section \ref{sec:epsc} based on $\epsilon$-circulant matrices.

\begin{theorem}\label{thm:new}
Let $\mathcal T(U)$ be the matrix in \eqref{U} and assume it is a
subgenerator.
For $x=(x_i)_{i=1,\ldots,n-1}\in\mathbb C^{n-1}$ define
  $H(x)=U_0+\sum_{i=1}^{n-1}x_iU_i$.   If $|x_i|\le 1$,
  $i=1,\ldots,n-1$, then $\|e^{sH(x)}\|_\infty\le 1$  for any  $s\ge 0$. Moreover,
$\|G(H(x),E)\|_\infty\le\|E\|_\infty$ for any $m\times m$ matrix $E$.
\end{theorem}
\pf
  Define $\alpha=\max_i(-(U_0)_{ii})$, $\tilde H=H((1,\ldots,1))$ and
  $B=\tilde H+\alpha I$. From the choice of $\alpha$ it follows that
  $B\ge 0$.  We have
\[
\|e^{s\tilde H}\|_\infty=\|e^{sB-s\alpha I}\|_\infty=e^{-s\alpha}\|e^{sB}\|_\infty\le
e^{-s\alpha}e^{\|sB\|_\infty}
\]
where the latter inequality holds in view of \cite[Theorem
10.10]{highambook}.  Since $B\ge 0$ we may write that
$\|B\|_\infty=\|B\uno\|_\infty$. From the inequality
$(\sum_{k=0}^{n-1}U_k)\uno\le 0$ we find that $(\alpha
I+\sum_{k=0}^{n-1}U_k)\uno\le \alpha\uno$, that is $B\uno\le
\alpha\uno$ whence $\|sB\|_\infty\le s\alpha$. Therefore we conclude
that $\|e^{s\tilde H}\|_\infty\le 1$ and the first claim is proved.
\\
Now, concerning
$H(x)$ we have
\begin{equation}\label{expu}
e^{sH(x)}=e^{sH(x)+s\alpha I -s\alpha I}=e^{-s\alpha}e^{sH(x)+s\alpha I}.
\end{equation}
Since $|x_k|\le 1$, $ U_0+\alpha I\ge 0$ and $U_k\ge 0$,
$k=1,\ldots,n-1$, we have
\[
|H(x)+\alpha I|\le 
|U_0+\alpha I| +| \sum_{k=1}^{n-1} x_k U_k|\le
U_0+\alpha I +\sum_{k=1}^{n-1}  U_k=B.
\]
By monotonicity of the infinity norm, we have
$\|H(x)+\alpha I\|_\infty \le \| B\|_\infty$, 
\[
\|e^{sH(x)+s\alpha I}\|_\infty \le e^{\| sH(x)+s\alpha I \|_\infty} \le e^{\| sB\|_\infty}\le e^{s\alpha},
\]
and,  from \eqref{expu},
$  %\[
\|e^{sH(x)}\|_\infty = e^{-s\alpha} \|e^{sH(x)+s\alpha I}\|_\infty
%\le e^{-s\alpha}e^{\|sH(x)+s\alpha I\|_\infty} 
\le 1.
$  %\]
From \eqref{frex}, we have 
\[
\|G(H(x),E)\|_\infty\le\|E\|_\infty
\int_0^1 \|e^{(1-s)H(x)}\|_\infty\|e^{sH(x)}\|_\infty ds\le \|E\|_\infty 
\]
and the last claim follows.
\qed

\section{Fast computations with Toeplitz and circulant matrices}\label{sec:fastcomp}
In this section we recall some basic properties of block-Toeplitz and
block-circulant matrices, useful for our computational
analysis. We refer the reader to Bini and Pan~\cite{bpbook} and
Bini {\it et al.}~\cite{blmbook} for more details.
Given a matrix $V\in\C^{m\times n}$, we denote by $V^T$ and by $V^H$
the transpose matrix and the transpose conjugate matrix of $V$,
respectively. The conjugate of a complex number $z$ is denoted by~$\overline z$.

Let $\cu$ be the imaginary unit such that $\cu^2=-1$ and
$\omega_n=\cos\frac{2\pi}n+\cu\sin\frac{2\pi}{n}$ be a primitive $n$th
root of the unity. We denote by $F=(\omega_n^{ij})_{i,j=0,n-1}$ the
Fourier matrix.  Recall that $F$ is nonsingular, that $F^{-1}=\frac 1n
F^H$ and that, given a vector $v\in\C^n$, the application $v\to u=Fv$
defines the inverse discrete Fourier transform (IDFT) of $v$. We
assume that $n$ is an integer power of 2, so that the vector $u$ can
be computed by means of the FFT algorithm in $\frac 32 n\log_2 n$
arithmetic operations (ops).  The application $u\to v=\frac 1n F^H u$
is called Discrete Fourier Transform (DFT) and the vector  $v$ can be computed in $\frac 32 n\log_2 n+ n$~ops.

Given the $m\times m$ matrices $V_i$, $i=0,\ldots,n-1$, we denote by
$V=(V_i)_{i=0,n-1}$ the block-(column) vector with block-entries  $V_i$, $i=0,\ldots,n-1$.
Finally, we define $\mathcal F=F\otimes I_m$, where $\otimes$ is
the Kronecker product and $I_m$ the identity matrix of order $m$. This
way, for a block-column vector $V$ the matrix $U=\mathcal F V$ can be
computed by means of $m^2$ IDFTs with $\frac 32nm^2\log_2n$
ops. Similarly, given the matrix $U$, the block-vector $V=\frac
1n\mathcal F^H U$ can be computed with $\frac 32nm^2\log_2n+nm^2$ ops.

\subsection{Block-circulant matrices}
For the results in this section we refer the reader to the book
\cite{bpbook} and to the references cited therein.

Given the block-vector $U=(U_i)_{i=0,n-1}$, with $m\times m$ blocks,
the $n\times n$ block-circulant matrix $\CC(U)=(C_{i,j})_{i,j=0,n-1}$ associated with $U$ is the matrix with block-entries
\[
C_{i,j}=U_{j-i\mathrm{~mod~} n}
\]
so that $[U_0,\ldots,U_{n-1}]$ coincides with the first block-row of
$\CC(U)$ and the entries of any other block-row are obtained by the entries
of the previous block-row by a cyclic permutation which moves the last block
entry to the first position and shifts the remaining block-entries one place to
the right. For instance, for $n=4$ one has
\[
\CC(U)=\left[\begin{array}{ccccc}
U_0&U_1&U_2&U_3\\
U_3&U_0&U_1&U_2\\
U_2&U_3&U_0&U_1\\
U_1&U_2&U_3&U_0\\
\end{array}\right].
\]
Observe that a block-circulant matrix is a particular block-Toeplitz matrix.

Block-circulant matrices can be simultaneously block-diagonalized by means of FFT, that is,
\[
\frac 1n \mathcal F^H\CC(U)\mathcal F=\diag(V_0,\ldots,V_{n-1}),\quad V=\mathcal FU.
\]
This property shows that block-circulant matrices are closed under
matrix multiplication, i.e., they form a matrix algebra, moreover the
product of a circulant matrix and a vector can be computed by means of
 Algorithm \ref{alg:1}.  This algorithm performs the computation
with $2m^2$ FFTs and $n$ matrix multiplications. Since $2m^3-m^2$ ops
are sufficient to multiply two $m\times m$ matrices, the overall cost
of Algorithm \ref{alg:1} is $3nm^2\log_2 n+nm^2+(2m^3-m^2)n$ ops.

\begin{algorithm}\caption{Product of a block-circulant matrix and a block-vector}\label{alg:1}

\SetKwInOut{Input}{Input}\SetKwInOut{Output}{Output}
\Input{Two block-vectors $X=(X_i)_{i=0,n-1},U=(U_i)_{i=0,n-1}$}
\Output{The block-vector $Y=\CC(U)X$, $Y=(Y_i)_{i=0,n-1}$}
\BlankLine
%Compute 
$V=\mathcal FU$

%Compute 
$Z=\mathcal FX$

%Compute 
$W_i=V_iZ_i$, $i=0,\ldots,n-1$,  $W=(W_i)_{i=0,n-1}$

%Compute 
$Y=\frac 1n \mathcal F^HW$

\end{algorithm}

If the input block-vectors are real then the vectors $V=\mathcal FU$ and 
$Z=\mathcal F X$ have
a special structure, that is, the components $V_0,Z_0$ and
$V_{n/2},Z_{n/2}$ are real while $V_i=\overline V_{n-i}$, $Z_i=\overline
Z_{n-i}$, for $i=1,\ldots,n/2-1$.
% where $\overline s$ denotes the complex conjugate of the complex number $s$. 
In this case, the number of matrix
multiplications at step 3 of Algorithm \ref{alg:1} is reduced to
$n/2$.

\begin{remark}\label{rem:cir}\rm
Observe that the product of two circulant matrices may be
computed by means of a product of a circulant matrix and a vector by
means of Algorithm~\ref{alg:1}. In fact, since the last column of the
block-circulant matrix $\CC(U)$ is the block-vector $\wh
U=(U_{n-i-1})_{i=0,n-1}$, if $\CC(Y)=\CC(U)\CC(X)$ then we find that $\wh
Y=\CC(U)\wh X$, where $\wh X=(X_{n-i-1})_{i=0,n-1}$, $\wh
Y=(Y_{n-i-1})_{i=0,n-1}$.
\end{remark}

\subsection{Block-triangular Toeplitz matrices}
We denote by $U$ the block-vector $(U_i)_{i=0,n-1}$ and by $\TT(U)$ the
block-upper triangular block-Toeplitz matrix whose first row is $[U_0,\ldots,U_{n-1}]$.
For $n=4$, for instance, 
\[
\TT(U)=\left[\begin{array}{cccc}
U_0&U_1&U_2&U_3\\
0&U_0&U_1&U_2\\
0&0&U_0&U_1\\
0&0&0&U_0
\end{array}
\right].
\]
Block-upper triangular block-Toeplitz matrices are closed under matrix multiplication.

Consider the vector $\wt U$ of $2n$ components obtained by filling the
vector $U$ with zero blocks, and the block-vector $V=(V_i)_{i=0,n-1}$ such that $V_0=0$ and
$V_i=U_{n-i}$ for $i=1,\ldots,n-1$.  Then the matrix $\CC(\wt U)$ can be
partitioned as follows
\begin{equation}\label{tc}
\CC(\wt U)=\left[\begin{array}{cc}
\TT(U)&\LL(V)\\
\LL(V)&\TT(U)
\end{array}\right],
\end{equation}
where $\LL(V)$ is the block-lower triangular block-Toeplitz matrix whose first block-column is $V$.
This expression enables one to compute the product $Y=\TT(U)X$ of a
block-upper triangular Toeplitz matrix and a block-vector with a low number of arithmetic operations. In fact, from \eqref{tc}
one deduces that $Y$ coincides with the first half of
the block-vector $\wt Y=\CC(\wt U)\wt X$ where
$\wt X$ is the block-vector of length $2n$ obtained by
filling $X$ with zeros.
This fact leads to Algorithm \ref{alg:2} for computing the product of a
block-triangular block-Toeplitz matrix and a block-vector.
The cost of this algorithm is  $6nm^2\log_2 (2n)+2nm^2+2(2m^3-m^2)n$ ops.

\begin{algorithm}\caption{Product of a block-triangular 
block-Toeplitz matrix and a block-vector}\label{alg:2}
\SetKwInOut{Input}{Input}\SetKwInOut{Output}{Output}

\Input{Two block-vectors $X=(X_i)_{i=0,n-1}$, $U=(U_i)_{i=0,n-1}$}
\Output{The block-vector $Y=\TT(U)X$, $Y=(Y_i)_{i=0,n-1}$}
\BlankLine

Set $\wt X=(\wt X_i)_{i=0,2n-1}$ with $\wt X_i=X_i$ for $i=0,\ldots,n-1$, $\wt X_i=0$ for $i=n,\ldots,2n-1$

Set $\wt U=(\wt U_i)$ with $\wt U_i=U_i$ for $i=0,\ldots,n-1$, $\wt U_i=0$ for $i=n,\ldots,2n-1$

Apply Algorithm \ref{alg:1} to compute $\wt Y=\CC(\wt U)\wt X$

Set $Y=(Y_i)_{i=0,n-1}$ with $Y_i=\wt Y_i$, for $i=0,\ldots,n-1$.

\end{algorithm}

\begin{remark}\label{rem:tt}\rm
Observe that the product of two block-triangular block-Toeplitz
matrices can be computed by means of a product of a block-triangular
block-Toeplitz matrix and a block-vector by means of Algorithm
\ref{alg:2}. In fact, since the last column of $\TT(U)$ is the block
vector $\wh U=(U_{n-i-1})_{i=0,n-1}$, if $\TT(Y)=\TT(U)\TT(X)$ then we find that $\wh
Y=\TT(U)\wh X$, where $\wh X=(X_{n-i-1})_{i=0,n-1}$, $\wh Y=(Y_{n-i-1})_{i=0,n-1}$.
\end{remark}

\subsection{Block-$\epsilon$-circulant matrices}
Given a block-vector $U$ and a complex number $\epsilon$, the block-$\epsilon$-circulant matrix $\CC_\epsilon(U)=(C_{i,j})$ is defined by 
\[
C_{i,j}=\left\{\begin{array}{ll}
U_{j-i}&\hbox{for } j\ge i,\\
\epsilon U_{n+j-i}&\hbox{for }j<i.
\end{array}\right.
\]
For instance, for $n=4$ one has
\[
\CC_\epsilon(U)=\left[\begin{array}{ccccc}
U_0&U_1&U_2&U_3\\
\epsilon U_3&U_0&U_1&U_2\\
\epsilon U_2&\epsilon U_3&U_0&U_1\\
\epsilon U_1&\epsilon U_2&\epsilon U_3&U_0\\
\end{array}\right].
\]
Observe that a block-$\epsilon$-circulant matrix is a particular case of
block-Toeplitz matrix and that, for $|\epsilon|$ small, a
block $\epsilon$-circulant matrix is an approximation of a block-triangular
block-Toeplitz matrix.

Like block-circulant matrices, block-$\epsilon$-circulant matrices
can be simultaneously block-diagonalized by means of FFT, so that they
are closed under matrix multiplication and form a matrix algebra as well. In
fact, one can show that
\begin{equation}\label{eq:epscirc}
\frac 1n \mathcal F\D_\epsilon^{-1}
\CC_\epsilon(U)\D_\epsilon\mathcal F^H=\diag(V_0,\ldots,V_{n-1}),\quad 
V=\mathcal F^H\D_\epsilon U,
\end{equation}
where
\[
\D_\epsilon=D_\epsilon\otimes I_m,~D_\epsilon=\diag(1,\theta ,\theta^2,\ldots,\theta^{n-1}),\quad
  \theta=\epsilon^{1/n}.
\]

\section{The exponential of a block-$\epsilon$-circulant matrix}\label{sec:expcirc}
Let $U=(U_i)_{i=0,n-1}$ be a block-vector of length $n$ where
$U_i\in\C^{m\times m}$, consider the block-$\epsilon$-circulant matrix
$\CC_\epsilon(U)$ and its matrix exponential
$e^{\CC_\epsilon(U)}$. In view of \eqref{eq:epscirc}, we find
that 
\[
e^{\CC_\epsilon(U)}=\frac 1n\D_\epsilon\mathcal F^H \diag(e^{V_0},\ldots,e^{V_{n-1}})\mathcal F \D_\epsilon^{-1}, \quad 
V=\mathcal F^H\D_\epsilon U.\\
\]
Therefore the exponential of a block-$\epsilon$-circulant matrix is still 
block-$\epsilon$-circulant. Moreover, we have $e^{\CC_\epsilon(U)}=\CC_\epsilon(Y)$ where
\begin{equation}\label{eq:expepsc}
Y=\frac 1n \D_\epsilon^{-1}\mathcal FW,\quad
W=(W_i)_{i=0,n-1},~W_i=e^{V_i},~i=0,\ldots,n-1,  
\end{equation}
and $V=\mathcal F^H\D_\epsilon U$.
The above equations allow to compute the exponential of an $n\times n$ block-$\epsilon$-circulant matrix by computing $n$ exponentials of $m\times m$ matrices and two Fourier transforms, as described in  Algorithm \ref{alg:expepsc}.

\begin{algorithm}\caption{Exponential of a block-$\epsilon$-circulant matrix}\label{alg:expepsc}

\SetKwInOut{Input}{Input}\SetKwInOut{Output}{Output}\SetKwInOut{Computation}{Computation}

\Input{A complex number $\epsilon$, the block-vector $U=(U_i)_{i=0,n-1}$
  defining the first block-row of the $\epsilon$-circulant matrix $\CC_\epsilon(U)$ }
\Output{The block-vector $Y=(Y_i)_{i=0,n-1}$ such that $\CC_\epsilon(Y)=e^{\CC_\epsilon(U)}$}

\BlankLine
%Compute 
$Z=\D_\epsilon U$

%Compute 
$V=\mathcal F^H Z$

%Compute 
$W_i=e^{V_i}$, $i=0,\ldots,n-1$, and set $W=(W_i)_{i=0,n-1}$

%Compute 
$R=\frac1n \mathcal FW$

%Compute 
$Y=\D_\epsilon^{-1}R$

\end{algorithm}

Observe that the multiplication of $U$ by the diagonal matrix $\D_\epsilon$
at step 1 reduces to scaling the blocks $U_i$ by the scalar
$\theta_i$. The multiplication by $\D_\epsilon^{-1}$ at step 5 performs
similarly.  Therefore the overall cost of the algorithm is given by
$3m^2n\log_2 n+3m^2n$ ops plus the cost of computing $n$ 
exponentials of $m\times m$ matrices.

For $\epsilon=1$ the block-$\epsilon$-circulant matrix turns to a block-circulant matrix and Algorithm \ref{alg:expepsc} takes the simpler form described in Algorithm \ref{alg:expcirc}. The computational cost in this case 
is reduced to $3m^2n\log_2 n+m^2n$ ops plus the cost of computing $n$ 
exponentials of $m\times m$ matrices.

\begin{algorithm}\caption{Exponential of a block-circulant matrix}\label{alg:expcirc}

\SetKwInOut{Input}{Input}\SetKwInOut{Output}{Output}\SetKwInOut{Computation}{Computation}

\Input{The block-vector $U=(U_i)_{i=0,n-1}$ defining the first block-row of $\CC(U)$}
\Output{The block-vector $Y=(Y_i)_{i=0,n-1}$ such that $\CC(Y)=e^{\CC(U)}$}

\BlankLine
%Compute 
$V=\mathcal F^H U$

%Compute 
$W_i=e^{V_i}$, $i=0,\ldots,n-1$, and set $W=(W_i)_{i=0,n-1}$

%Compute 
$Y=\frac1n \mathcal FW$

\end{algorithm}

\subsection{Numerical stability}
Let $U=(U_i)_{i=0,n-1}$ be the block-vector defining the first block
row of the subgenerator $\TT(U)$.  We analyze the error generated by
computing the exponential of the block-$\epsilon$-circulant matrix
$\CC_\epsilon(U)$ by means of Algorithm \ref{alg:expepsc} in floating
point arithmetic, where $\epsilon\in\C$ with $|\epsilon|<1$.

Here and hereafter  $\fl(\cdot)$ denotes the result
computed in floating point arithmetic of the expression between
parenthesis. The symbol $\doteq$ denotes equality up to lower order
terms, and similarly the symbol $\dotleq$ stands for inequality up to lower
order terms. The symbol $\mu$ denotes the machine precision.

We recall the following useful fact  (see \cite[page 71]{highambook})
\begin{equation}\label{eq:fl}
\fl(x y)=xy(1+\eta),\quad |\eta|\le \frac {2\sqrt 2}{1-2\mu}\mu=: \beta \mu,
 \quad \beta\doteq 2\sqrt 2,
\end{equation}
for $x,y\in\C$,  and we use the following properties involving norms, where $v\in\C^n$
\begin{equation}\label{eq:norm}
\begin{split}
&\|v\|_\infty\le\|v\|_2,\quad \|v\|_2\le\sqrt n\|v\|_\infty,\quad \|v\|_2\le\|v\|_1,\\
&\|Fv\|_\infty\le n\|v\|_\infty,\quad \|Fv\|_2\le \sqrt n\|v\|_2.
\end{split}
\end{equation}

In order to perform the error analysis of Algorithm \ref{alg:expepsc}, we recall
the following result concerning FFT (see \cite[page 453]{highambook}).

\begin{theorem}\label{th:fft}
Let $x$ be a vector of $n$ components, $n=2^q$, $q$ integer, $y=Fx$,
where $F=(\omega_n^{ij})_{i,j=0,n-1}$ is the Fourier matrix. Let $\tilde y$ be the
vector obtained in place of $y$ by applying the Cooley-Tukey FFT
algorithm in floating point arithmetic with precision $\mu$ where the
roots of the unity are approximated by floating point numbers up to
the error $\nu$. Then
\[
\frac{\|y-\tilde y\|_2}{\|y\|_2}\le \frac { q\eta}{1- q \eta},\quad \eta=\nu+(\sqrt 2+\nu)\frac{4\mu}{1-4\mu}.
\]
\end{theorem}

In particular, with $\nu=\mu$ and performing a first-order error
analysis where we consider only the part of the error which is linear
in $\mu$ we have
\begin{equation}\label{eq:fft}
\frac{\|y-\tilde y\|_2}{\|y\|_2}\le \gamma \mu q,\quad \gamma\doteq 4\sqrt 2+1.
\end{equation}
Observe that, since $F^H=\overline F$, we may replace $F$ with $F^H$
in the statement of Theorem \ref{th:fft}.

We split Algorithm \ref{alg:expepsc} into three parts.  The first
part consists in computing the entries of the matrices $V_k$ by means
of steps 1 and 2, the second part consists in computing the entries of
$W_k=e^{V_k}$ and the third part is formed by the remaining steps 4
and 5.  The first and third part can be viewed as the collection of
$m^2$ independent computations applied to the entry $(r,s)$ of the
generic block for $r,s=1,\ldots,m$.  More specifically, given the pair $(r,s)$,
denote
$u=(u_k)$, $z=(z_k)$, $v=(v_k)\in\mathbb C^n$ the vectors whose components are $(U_k)_{r,s}$,
$(Z_k)_{r,s}$, $(V_k)_{r,s}$, $k=0,\ldots,n-1$, respectively. The
computation of $u$ is obtained in the following way:
$\theta=\epsilon^{1/n}$, $z_k=\theta^k u_k$, for $k=0,\ldots,n-1$,
$v=F^Hz$.  While, denoting $w,r,y\in\C^n$ the vectors whose components
are $(W_k)_{r,s}$, $(R_k)_{r,s}$, $(Y_k)_{r,s}$, $k=0,\ldots,n-1$,
respectively, the computation of $y$ is obtained in the following way:
$r=\frac 1n F w$, $y_k=\theta^{-k}r_k$ for $k=0,\ldots,n-1$.

Define $\delta_z=\wt z- z$, $\delta_v=\wt v-v$, $\delta_r=\wt r-r$,
$\delta_y=\wt y-y$ where $\wt z,\wt v,\wt r, \wt y$ are the values
obtained in place of $z,v,r,y$ by performing computations in floating
point arithmetic. We denote also by $(\delta_r)_k=\wt r_k-r_k$ and
$(\delta_y)_k=\wt y_k-y_k$ the $k$-th component of $\delta_r$ and
$\delta_y$, respectively.

In our analysis we assume that the constants $\theta^k$ have been
precomputed and approximated with the numbers $\theta_k$ such that
$\theta_k= \theta^k(1+\sigma_k)$, $|\sigma_k|\le \mu$,
$k=-n+1,\ldots,0,\ldots,n-1$.

Since $z_k=\theta^ku_k$,
from \eqref{eq:fl}  we find that
 $\fl(\theta^ku_k)=\theta^ku_k(1+\eta_{k})(1+\sigma_{k})\doteq
 \theta^ku_k(1+\eta_{k}+\sigma_k)$. Thus,
\begin{equation}\label{eq:deltaz}
\|\delta_z\|_\infty\le  \mu\zeta\|z\|_\infty,\quad \zeta\doteq \beta+1. 
\end{equation}

Denoting by $\delta'$ the error introduced in computing the FFT $v$  of $z$ 
 in floating point arithmetic, we have  
\[
\delta_v=F^H \delta_z+\delta',
\]
and in view of  \eqref{eq:norm}, \eqref{eq:fft}, and \eqref{eq:deltaz} 
we obtain
\[\begin{split}
\|\delta_v\|_\infty&\le n\|\delta_z\|_\infty+\|\delta'\|_2
\le n\|\delta_z\|_\infty+\mu\gamma \log_2n \|v\|_2\\
&\le \mu\zeta n\|z\|_\infty +\mu \gamma   n\log_2n \|z\|_\infty\\
&=
\mu  n\|z\|_\infty(\zeta+\gamma \log_2n)\\
&\le  \mu n\|u\|_\infty(\zeta +\gamma \log_2n),
\end{split}\]
where the last  inequality follows from the fact that 
$\|z\|_\infty\le\|u\|_\infty$ since $z_k=\theta^k u_k$ and $|\theta|<1$.
This implies that $\Delta_{V_k}=\wt V_k-V_k$ is such that 
\[
\max_k |(\Delta_{V_k})_{r,s}|\le \mu n(\zeta+\gamma\log_2 n)\max_k |(U_k)_{r,s}|,
\]
which yields
\begin{equation}\label{eq:deltau}
\|\Delta_{V_k}\|_\infty\le m \max_k |(\Delta_{V_k})_{r,s}|\le\mu m n(\zeta+\gamma \log_2n)\max_{r,s,h}|(U_h)_{r,s}|.
\end{equation}

Concerning the second part of the computation, for the matrix $\Delta_{W_k}=\wt W_k-W_k$ we have
\begin{equation}\label{eq:deltaw}
\Delta_{W_k}=\fl(e^{\wt V_k})-e^{V_k}, ~~\fl(e^{\wt V_k})=e^{\wt V_k}+E_k
\end{equation}
where $E_k$ is the error generated by computing the matrix exponential
in floating point arithmetic.  Here we assume that $\|E_k\|_\infty\le
\mu\tau \|W_k\|_\infty$ for some positive constant $\tau$ which depends on the algorithm used to compute the matrix exponential.  
From the properties of the G\^{a}teaux derivative one has 
 $ \|e^{\wt V_k}-e^{V_k}\| \doteq \|G(V_k,\Delta_{V_k})\| $, and
from Theorem \ref{thm:new}, applied with $x_i=\bar\omega_n^{ik}
\theta^i$, $i=1,\ldots,n-1$, it follows that
$\|G(V_k,\Delta_{V_k})\|_\infty \le \|\Delta_{V_k}\|_\infty$
and $\|W_k\|_\infty\le 1$.

Combining these results with \eqref{eq:deltaw} leads to the bound
\begin{equation}\label{eq:deltaw1}
\|\Delta_{W_k}\|_\infty\dotleq \|\Delta_{V_k}\|_\infty+\mu\tau.
\end{equation}

Finally, for the third part of the computation, consisting of steps 4 and 5, we have
\[
\delta_r\doteq \frac 1n F \delta_w+\frac 1n \delta'' 
\]
where $\delta''$ is the error obtained by computing $Fw$ in floating point arithmetic. Thus from \eqref{eq:fft} we have
\begin{equation}\label{eq:deltar}
\begin{split}
\|\delta_r\|_\infty&\dotleq\frac 1n \|F\delta_w\|_\infty+ \mu \gamma \log_2n\|r\|_2\\
&\le \|\delta_w\|_\infty+\mu\gamma \sqrt n\log_2n\|y\|_\infty,
\end{split}
\end{equation}
where the second inequality holds from \eqref{eq:norm} and from $r_k=\theta^k y_k$ since $|\theta|\le 1$.

Moreover, %by setting $\delta_{y_k}=y_k-\widetilde y_k$, 
we find that
\begin{equation}\label{eq:deltay}
(\delta_y)_k=\theta^{-k}(\delta_r)_k+\theta^{-k}r_k\nu_k=
\theta^{-k}((\delta_r)_k+y_k\nu_k) ,\quad |\nu_k|\le\zeta\mu.
\end{equation}

Now we are ready to combine all the pieces and obtain the error bound on the computed value $Y$. From \eqref{eq:deltay} we get
\[
|(\delta_y)_k|\le
|\theta|^{-k}(\|\delta_{r}\|_\infty+\zeta \mu\|y\|_\infty).
\]
On the other hand, by using \eqref{eq:deltar}, we find that
\[
|(\delta_y)_k|\dotleq |\theta|^{-k}(\|\delta_w\|_\infty
  +(\zeta+\gamma \sqrt n\log_2n)\mu\|y\|_\infty).
\]
Thus we have
\[
\|\Delta_{Y_k}\|_\infty\le m|(\delta_y)_k|\dotleq
m|\theta|^{-k}(\max_h\|\Delta_{W_h}\|_\infty+
(\zeta+\gamma \sqrt n\log_2n)\mu
\max_h\|Y_h\|_\infty).
\]
Moreover, from \eqref{eq:deltaw1} and \eqref{eq:deltau} we conclude with the following bound
\[
\|\Delta_{Y_k}\|_\infty\dotleq m|\theta|^{-k}( \max_h\|\Delta_{V_h}\|_\infty+\mu\tau
+(\zeta+\gamma \sqrt n\log_2n)\mu
\max_h\|Y_h\|_\infty).
\]
Whence 
\[
\|\Delta_{Y_k}\|_\infty\dotleq \mu m|\theta|^{-k}( m n(\zeta+\gamma \log_2n)\max_{r,s,h}
|(U_h)_{r,s}|+\tau+ (\zeta+\gamma \sqrt n\log_2n)
\max_h\|Y_h\|_\infty)
\]
and we may conclude with the following

\begin{theorem}\label{th:err}
Let $\widehat Y_k$ be the value of $Y_k$ provided by Algorithm
\ref{alg:expepsc} applied in floating point arithmetic with precision
$\mu$ for computing $\CC_\epsilon(Y)=e^{\CC_\epsilon (U)}$, where
$Y=(Y_k)_{k=0,n-1}$, $U=(U_k)_{k=0,n-1}$.  Denote
$\Delta_{Y_k}=Y_k-\widehat Y_k$. One has
\[
\|\Delta_{Y_k}\|_\infty\dotleq \mu\epsilon^{-1} m\varphi
\]
where
\[
\varphi= m n(\zeta+\gamma \log_2n)\max_{r,s,h}
|(U_h)_{r,s}|+ (\zeta+\gamma \sqrt n\log_2n)
\max_h\|Y_h\|_\infty+\tau,
\]
$\zeta\doteq 1+2\sqrt 2$, $\gamma\doteq 4\sqrt 2 +1$, and $\tau\mu$ is the error bound in the computation of the matrix exponential, i.e., such that
$\|\fl(e^{V})-e^V\|_\infty \le \mu\tau\|e^V\|_\infty$ for an $m\times m$ matrix~$V$.
\end{theorem}

In the case where $\epsilon=1$,  we apply Algorithm \ref{alg:expcirc}
to compute the exponential of a block-circulant matrix and the above result leads to

\begin{theorem}\label{th:err2}
Let $\widehat Y_k$ be the value of $Y_k$ provided by Algorithm
\ref{alg:expcirc} applied in floating point arithmetic with precision
$\mu$ for computing $\CC(Y)=e^{\CC(U)}$, where
$Y=(Y_k)_{k=0,n-1}$, $U=(U_k)_{k=0,n-1}$.  Denote
$\Delta_{Y_k}=Y_k-\widehat Y_k$.   One has
\[
\|\Delta_{Y_k}\|_\infty\le \mu m\chi
\]
where $\chi=m n \gamma \log_2n \max_{r,s,h}
|(U_h)_{r,s}|+ \gamma \sqrt n\log_2n \max_h\|Y_h\|_\infty + \tau$, and
$\zeta\doteq 1+2\sqrt 2$, $\gamma\doteq 4\sqrt 2 +1$, 
and $\tau\mu$ is the error bound in the computation of the matrix exponential, i.e., such that
$\|\fl(e^{V})-e^V\|_\infty \le \mu\tau\|e^V\|_\infty$ for an $m\times m$ matrix $V$.
\end{theorem}

\section{The exponential of a block-triangular block-Toeplitz matrix}\label{sec:exptt}
Let $U=(U_i)_{i=0,n-1}$ be the block-vector defining the first block-row of the
subgenerator $\TT(U)$ of \eqref{U}.  
Since block-triangular block-Toeplitz matrices form a matrix algebra,
by using the Taylor series expansion of the matrix exponential, it follows that
$e^{\TT(U)}$ is still a block-triangular block-Toeplitz matrix. Denote by
$A=(A_i)_{i=0,n-1}$ the block-vector defining the entries on the first block-row of $e^{\TT(U)}$, i.e.,
such that $\TT(A)=e^{\TT(U)}$. In particular, we have $A_0=e^{U_0}$.

Let $K\ge n$ and define the $K$ dimensional block-vector
$U^{(K)}$ obtained by completing $U$ with zeros:
\begin{equation}\label{eq:uk}
U^{(K)}=(U_i^{(K)})_{i=0,K-1},\quad U_i^{(K)}=\left\{\begin{array}{ll}
U_i&\hbox { for }i=0,\ldots,n-1,\\
0 & \hbox { for }i=n,\ldots,K-1.
\end{array}\right.
\end{equation}
Consider the $K\times K$ block-triangular block-Toeplitz matrix
$\TT(U^{(K)})$. In view of \cite[Theorem 3.6]{highambook},
if $K_2>K_1\ge n$, then
$e^{\TT(U^{(K_1)})}$ is the principal $K_1\times K_1$ block-submatrix of
$e^{\TT(U^{(K_2)})}$. Denote
by 
$A^{(K)}=(A_i)_{i=0,K-1}$ the block-vector defining the first block-row of
$e^{\TT(U^{(K)})}$,
i.e., $A^{(K)}$ is the block-vector such that $\TT(A^{(K)})=e^{\TT(U^{(K)})}$. 

Let $\hat U=(\hat U_i)_{i=0,n-1}$ be such that 
\begin{equation}\label{uhat}
\hat U_0=U_0+\alpha I,~~\hat U_i=U_i,~~i=1,\ldots,n-1,
\end{equation}
where $\alpha=\max_j(-(U_0)_{j,j})$.
Define the block-vector
$\hat U^{(K)}$ with block-components $\hat U_i^{(K)}=\hat U_i$ for $i=0,\ldots,n-1$, and
$\hat U_i^{(K)}=0$ for $i=n,\ldots,K-1$.
 Observe that
$\mathcal T(\hat U^{(K)})=\mathcal T(U^{(K)})+\alpha I$ is a nonnegative matrix, and we may write
$e^{\mathcal T(U^{(K)})}=e^{-\alpha}e^{\TT(\hat U^{(K)})}$. We denote by 
$\hat A^{(K)}=(\hat A_i)_{i=0,K-1}$ the block-vector such that $\TT(\hat
A^{(K)})=e^{\TT(\hat U^{(K)})}$. In particular we have
$A_i=e^{-\alpha}\hat A_i$, $i=0,\ldots,K-1$.

\subsection{Decay properties}
In this section we investigate decay properties of the exponential
$e^{\TT(U^{(K)})}$ of a subgenerator, in the case where the
subgenerator is a banded block-triangular block-Toeplitz matrix. These
properties will be used in Section \ref{sec:emb} to estimate the
approximation error of the numerical method based on the embedding
into a block-circulant matrix. 

Decay properties of matrix functions have been analyzed in the
literature.  We refer the reader to the survey paper \cite{BeBoRa} and
to \cite{BeBo14}.  In our case the structure and sign properties play
an important role.  The matrix exponential $e^{\TT(U^{(K)})}$ is not
banded in general, but its off-diagonal entries have useful decay
properties for $K\to\infty$.  To prove this fact we need the following
result \cite[Theorem 3.6]{blmbook} on decay properties of analytic
functions.

\begin{theorem}\label{thm:dec}
  Let $H(z)=\sum_{i=0}^\infty z^i H_i$ be an $m\times m$ matrix power series
  analytic for $z\in\C$ with $|z|<R$, and $R>1$. For  any $1<\sigma<R$, the
  block-coefficients satisfy  
\begin{equation}\label{eq:dec}
|H_i|\le M(\sigma)\sigma^{-i},~~i=0,1,\ldots
\end{equation}
where $M(\sigma)$ is the $m\times m$ matrix with elements
$\max_{|z|=\sigma}|h_{rs}(z)|$, for $r,s=1,\ldots,m$, and the inequality \eqref{eq:dec} is meant componentwise. 
\end{theorem}

The following result provides bounds to $A_i$, $i=0,\ldots,K-1$.

\begin{theorem}\label{thm:decexp}
Let $K\ge n$ and let $\TT(A^{(K)})=e^{\TT( U^{(K)})}$, with $A^{(K)}=(
A_i)_{i=0,K-1}$, where $\mathcal T(U)$ in \eqref{U} is a subgenerator and $U^{(K)}$ is defined in \eqref{eq:uk}. For any $\sigma>1$, we have
\[
 A_i\uno\le e^{\alpha(\sigma^{n-1}-1)}\sigma^{-i}\uno,~~i=0,\ldots,K-1,
\]
where $\alpha=\max_j(-(U_0)_{j,j})$.
\end{theorem}

\pf
  We associate with the block-vector $\hat U=(\hat U_i)_{i=0,n-1}$ of
  \eqref{uhat} the $m\times m$ matrix polynomial $\hat
  U(z)=\sum_{h=0}^{n-1}z^h \hat U_h$. For the properties of block
  triangular block-Toeplitz matrices \cite{blmbook}, the matrix
  $\TT(\hat U^{(K)})^j$ is still a block-triangular block-Toeplitz
  matrix and the blocks in its first row are the coefficients of the
  matrix polynomial $P^{(j)}(z)=\hat U(z)^j~\mod~z^K$.  Let
  $P_i^{(j)}$ be the matrix coefficient of degree $i$ of $P^{(j)}(z)$,
  for $i=0,\ldots,K-1$.  From the power series expression of the
  matrix exponential we find that
\begin{equation}\label{eq:hvh}
\hat A_i=\sum_{j=0}^\infty \frac{1}{j!} P_i^{(j)},~~i=0,\ldots,K-1.
\end{equation}
We want to give an upper bound to the matrices $ P_i^{(j)}$.
Since $\hat U(z)^j$ is a matrix
polynomial, then it is analytic in all the complex plane and we may apply
Theorem \ref{thm:dec} with $H(z)=\hat U(z)^j$ and any $\sigma>1$. We have to
estimate the matrix $M(\sigma)$. The matrix coefficients of
$\hat U(z)$ are nonnegative, therefore for any $\sigma>1$ and for any $z\in\C$ with $|z|=\sigma$, we
have $|\hat U(z)^j|\le \hat U(\sigma)^j\le (\hat U(1)\sigma^{n-1})^j$.
Since $\mathcal T(U)$ is a subgenerator then $\hat U(1)\uno=(\alpha I+\sum_{h=0}^{n-1}U_h) \uno\le \alpha\uno$. So that we obtain 
$|\hat U(z)^j|\uno\le\alpha^j\sigma^{(n-1)j}\uno$.
Since $P^{(j)}(z)=\hat U(z)^j~\mod~z^K$, then $P_i^{(j)}$ is the matrix
coefficient of degree $i$ of $\hat U(z)^j$ and,
in view of \eqref{eq:dec}, we find that 
$P_i^{(j)}\uno\le\alpha^j\sigma^{(n-1)j}\sigma^{-i} \uno$.
From this inequality and from \eqref{eq:hvh} we obtain that for any
$\sigma>1$ and for $i=0,\ldots,K-1$
\[
\hat A_i\uno\le \sum_{j=0}^\infty \frac{1}{j!}\alpha^j \sigma^{(n-1)j} \sigma^{-i}\uno=\sigma^{-i}e^{\alpha\sigma^{n-1}}\uno.
\]
Since $\hat A_i=e^{-\alpha} A_i$ we conclude the proof.
\qed

\subsection{Method based on $\epsilon$-circulant matrix}\label{sec:epsc}
Let $\epsilon\in\mathbb C$ with $|\epsilon|$ sufficiently small,
and consider the block-$\epsilon$-circulant matrix
\begin{equation}\label{Ueps}
\CC_\epsilon(U)=\left[\begin{array}{cccc}
U_0 & U_1 & \ldots & U_{n-1}\\
\epsilon U_{n-1}    & U_0 & \ddots & \vdots \\
\vdots    &  \ddots   & \ddots & U_1 \\
 \epsilon U_{1}   &  \ldots   &   \epsilon U_{n-1}     & U_0
\end{array}\right].
\end{equation}

The exponential of $\CC_\epsilon(U)$ is still a block-$\epsilon$-circulant
matrix, that can be computed by means of Algorithm \ref{alg:expepsc}.  Denote
by $Y=(Y_i)_{i=0,n-1}$ the block-vector such that
$\CC_\epsilon(Y)=e^{\CC_\epsilon(U)}$.  The idea is to approximate the blocks
$A_i$, defining $\TT(A)=e^{\TT(U)}$, by the matrices $Y_i$, for
$i=0,\ldots,n-1$.

In order to estimate the approximation error, observe that 
the matrix $\CC_\epsilon(U)$ can be written as 
$\CC_\epsilon(U)=\TT(U)+\epsilon L$, where 
\begin{equation}\label{eq:L}
L=\left[\begin{array}{cccc}
0 & 0 & \ldots & 0\\
 U_{n-1}    & 0 & \ddots & \vdots \\
\vdots    &  \ddots   & \ddots & 0 \\
 U_{1}   &  \ldots   &  U_{n-1}     & 0
\end{array}\right].
\end{equation}
This property allows to give the following estimate:

\begin{theorem}\label{thm1}
Assume that $\TT(U)$ is a subgenerator, and that $\epsilon\in\C$.
One has
\[
\| e^{\mathcal T(U)}- e^{\CC_\epsilon(U)}\|_\infty \le 
e^{|\epsilon|  \|L\|_\infty}-1.
\]
Moreover, if $\epsilon$ is a pure imaginary number, then
\[
\| e^{\mathcal T(U)}-\mathrm{Re}(e^{\CC_\epsilon(U)})\|_\infty \le
e^{|\epsilon|^2  \|L\|_\infty^2}-1,
\]
where $\mathrm{Re}(e^{\CC_\epsilon(U)})$ is the real part of $e^{\CC_\epsilon(U)}$.
\end{theorem}

\pf
According to \eqref{tayl}, 
\begin{equation}\label{tayeps2}
e^{\CC_\epsilon(U)}=e^{\TT(U)}+\sum_{j=1}^\infty \frac{\epsilon^j}{j!}
G^{[j]}(\TT(U),L),
\end{equation}
where $G^{[j]}(\TT(U),L)$ are defined by means of 
\eqref{gatint}.
From  Proposition \ref{prop1} we obtain
\[
\| e^{\TT(U)}- e^{\CC_\epsilon(U)}\|_\infty  \le
\sum_{j=1}^\infty \frac{|\epsilon|^j}{j!}\| G^{[j]}(\TT(U),L)\|_\infty
\le \sum_{j=1}^\infty \frac{|\epsilon|^j}{j!} \|L\|_\infty^j=
e^{|\epsilon|\|L\|_\infty}-1.
\]
If $\epsilon$ is a pure imaginary number, since $e^{\TT(U)}$ is a real matrix, the
inequality is obtained by comparing the real parts in \eqref{tayeps2} and by
applying Proposition \ref{prop1}.
\qed

It is interesting to observe  that the choice of  an imaginary value for
$\epsilon$ provides an approximation error of the order
$O(|\epsilon|^2)$ instead of $O(|\epsilon|)$. The idea of using an
imaginary value for $\epsilon$ was used in \cite{high:complex} in the
framework of Fr\'echet derivative approximation of matrix functions.

The error bound can be improved by performing the computation with
several different values of $\epsilon$ and taking the mean of the real
parts of the results obtained this way. For instance, choose
$\epsilon_1=(1+\cu)\sqrt 2 \epsilon$, $\epsilon_2=-\epsilon_1$, where
$\epsilon>0$, and recall that $e^{\mathcal C_{\epsilon_j}(U)}$, $j=1,2$
are power series in $\epsilon$. Taking the arithmetic mean of
$e^{\mathcal C_{\epsilon_1}(U)}$ and $e^{\mathcal C_{\epsilon_2}(U)}$,
the components of odd degree in $\epsilon$ cancel out while the
coefficient of $\epsilon^2$ is pure imaginary. Therefore taking the real
part of the arithmetic mean provides an error $O(\epsilon^4)$. 

This technique can be generalized as follows. Choose an integer $k\ge
2$ and set $\epsilon_j=(\cu)^{1/k}\omega_k^j\epsilon$, $j=0,\ldots,k-1$, where 
$(\cu)^{1/k}$
is a principal $k$-th root of $\cu$. Then one can verify that the
arithmetic mean of $e^{\mathcal C_{\epsilon_j}(U)}$, $j=0,\ldots,k-1$ is
a power series in $\epsilon^k$, moreover, $\epsilon_j^k$ is a pure imaginary
number so that the real part of this mean provides an approximation
with error $O(\epsilon^{2k})$.

Observe that computing the exponential for different values of
$\epsilon$ might seem a substantial computational overload. However,
in a parallel model of computation, the exponentials $e^{\mathcal
  C_{\epsilon_j}(U)}$, $j=0,\ldots,k-1$, can be computed simultaneously
by different processors at the same cost of computing a single
exponential.

Algorithm \ref{alg:epscirc_int} reports this averaging technique.

\begin{algorithm}\caption{Exponential of a block-triangular matrix by means of $\epsilon$-circulant matrices and averaging}\label{alg:epscirc_int}

\SetKwInOut{Input}{Input}\SetKwInOut{Output}{Output}\SetKwInOut{Computation}{Computation}

\Input{The block-vector $U=(U_i)_{i=0,n-1}$ defining the first block-row of $\TT(U)$; a real number $\epsilon>0$; an integer $k>0$.}
\Output{The block-vector $Y=(Y_i)_{i=0,n-1}$ that is an approximation of the first block-row of
$e^{\TT(U)}$.}

\BlankLine

Set $\epsilon_j=(\cu)^{1/k}\omega_k^j\epsilon$, $j=0,\ldots,k-1$

%Compute 
Compute the first block row $W^{(j)}$ of $e^{\CC_{\epsilon_j}(U)}$,  $j=0,\ldots,k-1$, by means of Algorithm \ref{alg:expepsc}.

Set $Y=\frac1k \sum_{j=0}^{k-1}\mathrm{Re}(W^{(j)})$ 

\end{algorithm}

Theorem \ref{thm1} provides us with a bound on the error generated by
approximating the exponential of a block-upper triangular Toeplitz
matrix by means of the exponential of a block-$\epsilon$-circulant
matrix. In fact, in practical computations in floating point
arithmetic, the overall error is formed by two components: one
component is given by the approximation error analyzed in Theorem
\ref{thm1}, the second component is due to the roundoff and is
estimated by Theorem \ref{th:err}. More precisely, the effectively
computed approximation in floating point arithmetic is the
block-vector with components $\widehat Y_k=Y_k+\Delta_{Y_k}$,
$k=0,\ldots,n-1$, where $\|\Delta_{Y_k}\|_\infty$ is bounded in
Theorem \ref{th:err}.  On the other hand, $Y_k=A_k+E'_k $ where, by
Theorem \ref{thm1}, $E'_k$ is such that
\[
\|[E'_0,\ldots,E'_{n-1}]\|_\infty\le \left\{\begin{array}{ll}
\psi(|\epsilon|^2\| L\|_\infty^2)& \hbox{if~} 
\epsilon \hbox{~is imaginary}\\[1ex]
\psi(|\epsilon|\| L\|_\infty)&    \hbox{otherwise}
\end{array}\right.
\]
where $\psi(t)=e^t-1$.
This way, for the overall error $E_k=\Delta_{Y_k}+E_k'$ one has
\[
\|E_k\|_\infty\le\|\Delta_{Y_k}\|_\infty+\|E_k'\|_\infty\le m\mu\epsilon^{-1}\varphi+\psi(t),
\]
for $t=|\epsilon|\|L\|_\infty$, or $t=|\epsilon|^2\|L\|_\infty^2$.

This shows the need to find a proper balance between the two errors:
small values for $|\epsilon|$ provide a small approximation error
$\|E'_k\|_\infty$ but the roundoff errors diverge
to infinity as $\epsilon\to 0$. A good compromise is to choose $\epsilon$ so that the
upper bounds to $\|E'_k\|$ and $\|\Delta_{Y_k}\|_\infty$ have the same
order of magnitude.  Equating these upper bounds in the case of
non-imaginary $\epsilon$ yields
\[
|\epsilon|=\sqrt{m\mu\varphi/\|L\|_\infty},\quad \|E_k\|_\infty\dotleq 2\epsilon
\|L\|_\infty 
\]
and in the case of imaginary $\epsilon$,
\[
|\epsilon|=\sqrt[3]{m\mu\varphi/\|L\|_\infty^2},\quad \|E_k\|_\infty\dotleq
2\epsilon^2\|L\|_\infty^2. 
\]
The latter bound is an $O(\mu^{2/3})$. This implies that
asymptotically, as $\mu\to 0$, we may loose $1/3$ of the digits
provided by the floating point arithmetic.

If we adopt the strategy of performing the computation with $k$
different values of $\epsilon_j=(\cu)^{1/k}\omega_k^j\epsilon$,
$j=0,\ldots,k-1$, so that the approximation error is $O(\epsilon^{2k})$,
then the total error turns to $O(\mu^{2k/{2k+1}})$, i.e., only
$1/(2k+1)$ digits are lost.

An interesting point is that the quantities $\| L\|_\infty$ and
$\max_{r,s,h}|(U_h)_{r,s}|$ are involved in the expressions of the
error bound. Since $\TT(U)$ is a generator, both these quantities are
bounded from above by $\alpha=\max_j(-(U_0)_{j,j})$. However, by means
of simple manipulations, we may scale the input so that it is bounded
by 1. This is performed by applying to $\TT(U)$ the scaling and
squaring technique of \cite{higham09}.

Let $p\ge 0$ be an integer such that %$\| \TT(U)\|_\infty\le 2^p$. 
$\alpha\le 2^p$. 
Then, since $e^{\TT(U)}=(e^{\TT(U/2^p)})^{2^p}$, we first
compute $ e^{\TT(U/2^p)}$ and then recover $e^{\TT(U)}$ by performing
$p$ repeated matrix squaring.  
In this way we have $\| L/2^p\|_\infty<1$ and $\max_{r,s,h}|(U_h)_{r,s}/2^p|<1$.
Since $\TT(U/2^p)$ is still a generator, the error analysis performed for
$e^{\TT(U)}$ applies as well, and we can
approximate the first block-row of $ e^{\TT(U/2^p)}$
with the first block-row $Y=(Y_i)$ of $ e^{\TT(U_\epsilon/2^p)}$ for a suitable
$\epsilon\in\C$ with $|\epsilon|<1$.
Finally we recover an approximation to $e^{\TT(U)}$ by computing
$\TT(Y)^{2^p}$ by means of $p$ repeated squarings, by using the Toeplitz
structure and Algorithm \ref{alg:2}, in view of Remark \ref{rem:tt}.
The overall procedure is described in  Algorithm \ref{alg:ttepsc}.

\begin{algorithm}\caption{Exponential of a block-triangular block-Toeplitz matrix by using $\epsilon$-circulant matrices}\label{alg:ttepsc}

\SetKwInOut{Input}{Input}\SetKwInOut{Output}{Output}\SetKwInOut{Computation}{Computation}

\Input{The block-vector $U=(U_i)_{i=0,n-1}$ defining the first block-row of
  $\TT(U)$, $\epsilon\in\C$} 

\Output{The block-vector $Y=(Y_i)_{i=0,n-1}$,
  that is an approximation of the first block-row of $e^{\TT(U)}$}

\BlankLine

$\alpha=\max_j(-(U_0)_{j,j})$, $p=\lfloor \log_2\alpha\rfloor +1$ and 
$\tilde U=U/2^p$

Compute $Y$, the first block-row of $e^{\CC_\epsilon(\tilde U)}$, by means of Algorithm \ref{alg:expepsc} 

If $\epsilon$ is imaginary, replace $Y$ with the real part of $Y$

\For{$r=1,\ldots,p$}{compute  $\TT(Y)=\TT(Y)\TT(Y)$}

\end{algorithm}

\subsection{Embedding into a circulant matrix}\label{sec:emb}

The idea of this method is to embed the matrix $\mathcal T(U)$ into a $K\times
K$ block-circulant matrix $\mathcal C(U^{(K)})$. The first block-row of $e^{
  \mathcal T(U)}$ is approximated by the first $n$ blocks of the first block-row  of $e^{\mathcal C(U^{(K)})}$. 
Specifically, take $K\ge n$ and consider the block-vector $U^{(K)}$
defined in \eqref{eq:uk}. 
The block-circulant matrix $\mathcal C(U^{(K)})$ may be
partitioned as
\[
\mathcal C(U^{(K)})=\left[\begin{array}{cc}
\mathcal T(U) & P \\
Q & \mathcal T(U^{(K-n)})
\end{array}\right],
\]
where $P$ and $Q$ are $n\times (K-n)$ and $(K-n)\times n$ block-matrices,
respectively. 

Denote by $\mathcal E_1$ and by $\mathcal E_K$ the $(mnK)\times (mn)$ matrices formed by the first
$mn$ and the last $mn$ columns, respectively, of the identity matrix of size
$mnK$. 
The matrix $\mathcal C(U^{(K)})$ can be also written as
\begin{equation}\label{dec}
\mathcal C(U^{(K)})=\mathcal  T(U^{(K)})+H_K,~~~H_K=\mathcal E_K L \mathcal E_1^T,
%%BD controllare T_k
\end{equation}
where the matrix $L$ is defined in \eqref{eq:L}.  Because of the
triangular Toeplitz structure, the desired matrix $e^{\mathcal T(U)}$
is identical to the $n\times n$ block-leading submatrix of
$e^{\mathcal T_K(U)}$.
Our idea is to approximate the first block-row of
$e^{\mathcal T(U)}$ with the first $n$ blocks of the first row of
$e^{\mathcal C(U^{(K)})}$.  As pointed out in Section \ref{sec:expcirc},
$e^{\mathcal C(U^{(K)})}$ is a block-circulant matrix, and can
be computed by means of Algorithm \ref{alg:expcirc} with $3m^2K\log_2
K+m^2K$ ops, plus the cost of computing $K$ exponentials of $m\times
m$ matrices. 

Denote by
$S^{(K)}=(S_i^{(K)})_{i=0,K-1}$ the first block-row of $e^{\mathcal
  C(U^{(K)})}$, so that $\mathcal C(S^{(K)})=e^{\mathcal C(U^{(K)})}$.
An approximation of the matrices $A_i$, $i=0,\ldots,n-1$, defining the first block-row of $e^{ \mathcal T(U)}$
is provided by $S_i^{(K)}$, $i=0,\ldots,n-1$;
as $K$ increases, the approximation
improves, as shown by the following result.

\begin{theorem}\label{thm:erremb}
Let 
$e^{ \mathcal T(U)}= \mathcal T(A)$, with $A=(A_i)_{i=0,n-1}$. Let $K\ge n$
and let $e^{\mathcal C(U^{(K)})}=\mathcal C(S^{(K)})$, with
$S^{(K)}=(S_i^{(K)})_{i=0,K-1}$. 
One has $S_i^{(K)}-A_i\geq 0$ for $i=0,\ldots,n-1$,  and
\begin{equation}\label{boundemb}
\left\|  \left[\begin{array}{ccc}S_0^{(K)}-A_0 & \ldots &
    S_{n-1}^{(K)}-A_{n-1}\end{array}
\right] \right\|_\infty \le  f_K(\sigma)
\end{equation}
for any $\sigma>1$, where
\[
f_K(\sigma)=
(e^{\|L\|_\infty} -1)e^{\alpha(\sigma^{n-1}-1)}\frac{\sigma^{-K+n}}{1-\sigma^{-1}},
\]
with $\alpha=\max_j(-(U_0)_{j,j})$ and $L$ defined in \eqref{eq:L}.
\end{theorem}

\pf
By using \eqref{dec} and \eqref{tayl}, we find that 
%%BD T_k
\[
e^{\mathcal C(U^{(K)})}-e^{\mathcal T(U^{(K)})}=
\sum_{j=1}^\infty \frac{1}{j!} G^{[j]}(\mathcal  T(U^{(K)}),H_K).
\]
Equating the first $n$ blocks in the first block-row in the above equation yields
\begin{equation}\label{e1te1}
\left[\begin{array}{ccc}S_0^{(K)} -A_0 & \ldots &
    S_{n-1}^{(K)}-A_{n-1}\end{array}
\right]= %I_{m,mnK} ( \mathcal C(S^{(K)}) - \mathcal T(A))I_{mnK,k}=\\
\sum_{j=1}^\infty \frac{1}{j!} W^{[j]},
\end{equation}
where $W^{[j]}$ is the block-row vector formed 
by the first $n$ block-entries in the first block-row of 
$ G^{[j]}(\mathcal T_K(U),H_K)$. That is,
%%BD T_k
\[
W^{[j]}=\widehat{\mathcal E}_1^T  G^{[j]}(\mathcal  T(U^{(K)}),H_K)\mathcal E_1,
\]
where $\widehat{\mathcal E}_1$ is the $mnK\times m$ matrix formed by the first $m$ columns of the identity matrix.
Since $H_K\ge 0$, from \eqref{e1te1} and from Proposition \ref{prop1} we deduce that $W^{[j]}\ge 0$ so that
$S_i^{(K)}-A_i\ge 0$ for $i=0,\ldots,n-1$.
On the other hand, in view of \eqref{gatint} and from the fact that $H_K=\mathcal E_K L \mathcal E_1^T$, we may write
\begin{equation}\label{wj}
\begin{split}
%%BD
W^{[j]}  &= 
j \int_0^1 \widehat{\mathcal E}_1^T  
 e^{ (1-s)\mathcal  T(U^{(K)})} \mathcal E_K L \mathcal E_1^T
G^{[j-1]}(s\mathcal  T(U^{(K)}),H_K)  \mathcal E_1 ds\\
& =  j \int_0^1
V(s) L Z^{[j-1]}(s)ds
\end{split}
\end{equation}
%%BD T_K
where $V(s)=\left[\begin{array}{ccc}V_0(s) & \ldots &
    V_{n-1}(s)\end{array}\right]$ is the block-row vector formed by the last $n$ block-entries of 
the first block-row of
$ e^{ (1-s)\mathcal  T(U^{(K)})}$, and $ Z^{[j-1]}(s)$ is the $n\times n$ block
leading submatrix of  $G^{[j-1]}(s \TT(U^{(K)}),H_K)$.

Since the matrix $(1-s)\mathcal  T(U^{(K)})$ is a
subgenerator,   it follows
that $V_i(s)\ge 0$ and, from Theorem \ref{thm:decexp}, that for any $\sigma>1$,
\[
V_i(s)\uno \le e^{(1-s)\alpha(\sigma^{n-1}-1)}\sigma^{-K+n-i}\uno
\le e^{\alpha(\sigma^{n-1}-1)}\sigma^{-K+n-i}\uno,
\]
for $i=0,\ldots,n-1$, where the latter inequality follows from the fact that
$1-s\le 1$. This implies that
\[
\|V(s)\|_\infty\le e^{\alpha(\sigma^{n-1}-1)}\sum_{i=0}^{n-1}\sigma^{-K+n-i}\le 
e^{\alpha(\sigma^{n-1}-1)}\frac {\sigma^{-K+n}}{1-\sigma^{-1}}.
\]
%%BD T_K
 Moreover, since $s\TT(U^{(K)}) $ is a subgenerator, $H_K$ is
nonnegative and $\|H_K\|_\infty=\| L\|_\infty$, then, from Proposition
\ref{prop1}, we have $G^{[j-1]}(s\TT(U^{(K)}),H_K)\ge 0$ and
\[
\|G^{[j-1]}(s\TT(U^{(K)}),H_K)\|_\infty\le s^{j-1}\| L\|_\infty^{j-1}.
\]
This
latter inequality implies that $\| Z^{[j-1]}(s)\|_\infty\le s^{j-1}\|
L\|_\infty^{j-1}$.  Therefore, by taking norms in \eqref{wj}, we find that
\[
\begin{split}
%%BD Verificare 1/j e j
\| W^{[j]}\|_\infty & \le j
\int_0^1 \| V(s)\|_\infty \| L \|_\infty \| Z^{[j-1]}(s)\|_\infty ds \\
& \le
j \| L\|_\infty^{j} e^{\alpha(\sigma^{n-1}-1)}
\frac{\sigma^{-K+n}}{1-\sigma^{-1}}
%%BD verif = o <=
 \int_0^1 s^{j-1} ds = \| L\|_\infty^{j} e^{\alpha(\sigma^{n-1}-1)}\frac{\sigma^{-K+n}}{1-\sigma^{-1}}.
\end{split}
\]
Hence, by taking norms in \eqref{e1te1}, we obtain \eqref{boundemb}.
\qed

\begin{remark}The matrices $A_i$ and $S_i^{(K)}$ have a probabilistic interpretation. Namely, the matrix $A_i$ is the
probability that the BMAP is absorbed after time 1, {\em and} at time
1 there have been $i < n$ arrivals; the 
 matrix $S_i^{(K)}$ is the probability that the BMAP is absorbed
after time 1, {\em and} at time 1 there have been $i$, or $i+K$, or
$i+2 K$, or \ldots, arrivals.   
Clearly, there are more trajectories favourable for $S_i^{(K)}$ than
for $A_i$ and $A_i \leq  S_i^{(K)}$.    Similarly, there are
more trajectories favourable for $S_i^{(K)}$ than for $S_i^{(\ell
  K)}$ for a positive integer $\ell$.  This shows that,  if we take a
sequence of integers $\ell_1$, $\ell_2$, \ldots, and a sequence $K_0$,
$K_1$, $K_2$, \ldots, such that $K_{n+1} = \ell_{n+1} K_n$, then
\[
S_i^{(K_0)} \geq S_i^{(K_1)} \geq S_i^{(K_2)} \geq \cdots \geq A_i
\]
for $i=0,\ldots,K_0-1$. Therefore, the sequence $\{S^{(K)}\}$ has some 
monotonicity  property in its convergence to $A$.
\end{remark}

The bound in \eqref{boundemb} shows that the error has an exponential
decay as $K$ increases. Moreover, such bound holds for any $\sigma
>1$. Therefore we can fix a tolerance $\epsilon$ and a $\sigma>1$, and
find $K$ such that $f_K(\sigma)<\epsilon$. Since we would like to keep
$K$ as low as possible, another way to proceed is to fix a tolerance
$\epsilon$ and find $\sigma$ such that the size $K$ for which
$f_K(\sigma)<\epsilon$ is minimum. More specifically, after some
manipulations, from the condition $f_K(\sigma)<\epsilon$ we obtain that
$K>g(\sigma)$ where
\[
g(\sigma)=\frac{
\alpha(\sigma^{n-1}-1)+\log(\sigma/(\sigma-1))+ 
\log(\epsilon^{-1})+\log(e^{\| L\|_\infty}-1)
}
{\log(\sigma)}+n.
\]
Since $\sigma>1$ is arbitrary, we choose $\sigma$ such that
$g(\sigma)$ has a minimum value. In fact, the function $g(\sigma)$
diverges to infinity as $\sigma$ tends to 1 and to $\infty$, therefore
it has at least a local minimum $\sigma^*$ and we can choose $K>g(\sigma^*)$.

When we perform the computation in floating point arithmetic, we have to
consider also the error generated by roundoff in computing the exponential of
a block-circulant matrix. In practical computations, we obtain
a block-vector with components $\widehat
Y_i=Y_i+\Delta_{Y_i}$, $k=0,\ldots,n-1$, where
$\|\Delta_{Y_i}\|_\infty$ is bounded in Theorem \ref{th:err2} and $Y_i=A_i+E'_i $ where, by Theorem \ref{thm:erremb}, $E'_i$ is
such that
\[
\|[E'_0,\ldots,E'_{n-1}]\|_\infty\le  f_K(\sigma).
\]
Altogether, for the overall error $E_i=\Delta_{Y_i}+E_i'$, one has
\[
\|E_i\|_\infty\le\|\Delta_{Y_i}\|_\infty+\|E_i'\|_\infty\le m\mu \chi + f_K(\sigma).
\]

A similar analysis can be carried out for the relative error. In this
case the inequality $f_K(\sigma)<\epsilon$ is replaced by
$f_K(\sigma)<\hat\epsilon$, for $\hat\epsilon=\epsilon\| [A_0,\ldots,A_{n-1}] \|_\infty$. So that the function $g(\sigma)$ is modified by replacing $\epsilon$ with $\hat\epsilon$.

Like at the end of Section \ref{sec:epsc}, in the overall estimate of
the error, the quantities $\| L\|_\infty$ and
$\max_{r,s,h}|(U_h)_{r,s}|$ are bounded from above by
$\alpha=\max_j(-(U_0)_{j,j})$, and we may scale the block-vector $U$
so that these quantities are bounded by 1. 
% In fact, let $p\ge 0$ be an
% integer such that %$\| \TT(U)\|_\infty\le 2^p$.
% $\alpha\le 2^p$.  Therefore, we first approximate $ e^{\TT(U/2^p)}$
% and then recover $e^{\TT(U)}$ by performing $p$ repeated matrix
% squarings.

The overall procedure is summarized in Algorithm \ref{alg:ttemb},
where the repeated squaring of the block-triangular block-Toeplitz
matrices can be performed by using Algorithm~\ref{alg:2}, as explained
in Remark \ref{rem:tt}.

\begin{algorithm}\caption{Exponential of a block-triangular block-Toeplitz matrix by using embedding into a circulant matrix}\label{alg:ttemb}

\SetKwInOut{Input}{Input}\SetKwInOut{Output}{Output}\SetKwInOut{Computation}{Computation}

\Input{The block-vector $U=(U_i)_{i=0,n-1}$, an integer $K>n$}
\Output{The block-vector $Y=(Y_i)_{i=0,n-1}$, that is an approximation of the
first block-row of $e^{\TT(U)}$}

\BlankLine

Set $\alpha=\max_j(-(U_0)_{j,j})$, $p=\lfloor \log_2\alpha\rfloor +1$ and 
$\tilde U=U/2^p$

Set $W=(W_i)_{i=0,K-1}$ with $W_i=\tilde U_i$ for $i=0,\ldots,n-1$, $W_i=0$ for $i=n,\ldots,K-1$

Apply Algorithm \ref{alg:expcirc} to compute the first block-row $V=(V_i)_{i=0,K-1}$ of $e^{\CC(W)}$

Set $Y_i=V_i$, for $i=0,\ldots,n-1$.

\For{$r=1,\ldots,p$}{compute  $\TT(Y)=\TT(Y)\TT(Y)$}

\end{algorithm}

\subsection{Taylor series method}\label{sec:tay}
In this section we use the Taylor series method for computing the
exponential of an essentially nonnegative matrix, where the
block-triangular block-Toeplitz structure is exploited to perform fast
matrix-vector multiplications.  The computation of the exponential of
an essentially nonnegative matrix have been analyzed in \cite{xy13}
and \cite{shao}.

Following \cite{xy13} and \cite{shao}, the Taylor series method is
applied to compute $e^{\TT(\hat U)}$, since the matrix $\TT(\hat
U)=\TT(U)+\alpha I$ is nonnegative and $e^{\TT(U)}$ can be obtained by
means of the equation $e^{\TT(U)}=e^{-\alpha}e^{\TT(\hat U)}$. In this
way, we avoid possible cancellations in the Taylor summation.

Denote by $S_r(\TT(\hat U))$ the Taylor series truncated at the $r$th 
term, namely
\[
S_r(\TT(\hat U))=\sum_{k=0}^{r-1}\frac{\TT(\hat U)^k}{k!}.
\]
The following bound on the approximation error is given in \cite{xy13}. 

\begin{theorem}\label{thm:trtay}
Let $r$ be such that $\rho(\TT(\hat U)/(r+1))<1$. Then
\[
| e^{\TT(\hat U)}-S_r(\TT(\hat U))|\le \frac{\TT(\hat U)^r}{r!}\left(
I-\frac{\TT(\hat U)}{r+1}\right)^{-1}.
\]
\end{theorem}

%In particular, since $\left(
%I-\frac{\TT(\hat U)}{r+1}\right)^{-1}\ge I$, the condition
%$\frac{\TT(\hat U)^r}{r!}\le \epsilon S_r(\hat U)$ implies the following 
%bound on the
%approximation error
%$| e^{\TT(\hat U)}-S_r(\TT(\hat U))|\le \epsilon  S_r(\hat U)$.
%%
%\remarque{\rouge{Remarque from Guy: I do not understand this sentence}
%}
%%BD verificare

The scaling and squaring method is used to accelerate the convergence of the
Taylor series, by using the property that
\[
e^{\TT(U)}=e^{-\alpha}e^{\TT(\hat U)}=\left(e^{-\alpha/2^p}e^{\TT(\hat
    U)/2^p}\right)^{2^p}.
\]
Indeed, if $\tilde \rho$ is an estimate of $\rho(\TT(\hat U))$, and if
$p=\lfloor \log_2\tilde\rho\rfloor +1$, then $\rho(\TT(\hat U)/2^p)<1$ and the
truncated Taylor series expansion is used to approximate $e^{\TT(\hat
  U)/2^p}$. 
Since $\TT(\hat U)$ is block-triangular block-Toeplitz, then
$\rho(\TT(\hat U))=\rho(\hat U_0)$.

The Toeplitz structure is used in the computation of the Taylor expansion and
in the squaring procedure. In fact, the computation of each term  in the power series expansion consists in performing products
between block-triangular block-Toeplitz matrices, that can be done by
applying Algorithm \ref{alg:2} in view of Remark \ref{rem:tt}; similarly in
the squaring procedure at the end of the algorithm.

Concerning rounding errors, we observe that the Taylor polynomial is
the sum of nonnegative terms. Therefore no cancellation error is
encountered in this summation. The main source of rounding errors is
the computation of the powers $\TT(\hat U)^k$ for $k=2,\ldots,$ which
are computed by means of Algorithm \ref{alg:2} in view of Remark
\ref{rem:tt} relying on FFT.  We omit the error analysis of this
computation, which is standard. However, we recall that in view of
Theorem \ref{th:fft}, FFT is normwise backward stable but not
component-wise stable. For this reason, for the truncation of the
power series it is convenient to replace the component-wise bound
expressed by Theorem \ref{thm:trtay} by the norm-wise bound
\[
\| e^{\TT(\hat U)}-S_r(\TT(\hat U))\|_\infty\le 
\left\|\frac{\TT(\hat U)^r}{r!}\left(
I-\frac{\TT(\hat U)}{r+1}\right)^{-1}\right\|_\infty,
\]
from which we obtain that the condition
$
\left\|\frac{\TT(\hat U)^r}{r!}\left(
I-\frac{\TT(\hat U)}{r+1}\right)^{-1}\right\|_\infty < 
\epsilon \left\|S_r(\TT(\hat U))\right\|_\infty
$
implies that $\| e^{\TT(\hat U)}-S_r(\TT(\hat U))\|_\infty\le 
\epsilon \|S_r(\TT(\hat U))\|_\infty$.

The overall procedure is stated in
Algorithm \ref{alg:ttt}.

It is worth pointing out that, if the computation of the powers of the
triangular Toeplitz matrices is performed with the standard algorithm
then the computation is component-wise stable as shown in \cite{xy13}.

\begin{algorithm}\caption{Exponential of a block-triangular block-Toeplitz
    matrix by using Taylor series expansion}\label{alg:ttt}

\SetKwInOut{Input}{Input}\SetKwInOut{Output}{Output}
\SetKwInOut{Computation}{Computation}

\Input{The block-vector $U=(U_i)_{i=0,n-1}$ defining the first block-row of 
$\TT(U)$, a tolerance $\epsilon>0$, a maximum
number of iterations $K$}
\Output{The block-vector $Y=(Y_i)_{i=0,n-1}$, that is an approximation of the
first block-row of $e^{\TT(U)}$}

\BlankLine

Set $\alpha=\max_j(-(U_0)_{j,j})$

Set $\hat U=(U_i)_{i=0,n-1}$, $\hat U_0=U_0+\alpha I$, $\hat U_i=U_i$,
$i=1,\ldots,n-1$

Compute $\tilde \rho$ an estimate of $\rho(\hat U_0)$, or set $\tilde
\rho=\|\hat U_0\|_\infty$

Compute $p=\lfloor \log_2\tilde\rho\rfloor +1$ and $V=\hat U/2^p$

Set $W=V$ and $Y=(Y_i)_{i=0,n-1}$, $Y_0=I+V_0$, $Y_i=V_i$, $i=1,\ldots,n-1$.

\For{$r=2,\ldots,K$}{

Compute $\TT(W)=\TT(V)\TT(W/r)$ and $Y=Y+W$

\If{$\| W\|_\infty < \epsilon \| Y\|_\infty$}{break}
}
Compute $Y=e^{-\alpha/2^p}Y$

\For{$i=1,\ldots, p$}{compute $\TT(Y)=\TT(Y)\TT(Y)$}

\end{algorithm}

\section{Numerical experiments}\label{sec:exper}
The numerical experiments have been performed in Matlab. To compute the error
obtained with the proposed algorithms we have first computed the
exponential by using the \texttt{vpa} arithmetic of the Symbolic Toolbox
with 40 digits and we have considered this approximation as the exact
value.

Denote by $\tilde A_h$, $h=0,\ldots,n-1$, the approximations of the blocks on the first row of $e^{\TT(U)}$ and define the four errors 
\[
\begin{split}
&\textrm{cw-abs}=\max_{h,i,j}|(A_h)_{i,j}-(\tilde A_h)_{i,j}|,\\
&\textrm{cw-rel}=\max_{h,i,j}\{|(A_h)_{i,j}-(\tilde A_h)_{i,j}|/|(A_h)_{i,j}|\},\\
&\textrm{nw-abs}=\|[A_0-\tilde A_0,\ldots,A_{n-1}-\tilde A_{n-1}]\|_\infty,\\
&\textrm{nw-rel}=\|[A_0-\tilde A_0,\ldots,A_{n-1}-\tilde A_{n-1}]\|_\infty/\|[A_0,\ldots,A_{n-1}]\|_\infty,\\
\end{split}
\]
which represent absolute/relative component-wise and
norm-wise errors, respectively.

We compare the accuracy and the execution times of the proposed algorithms.

The test matrix $\TT(U)$ is taken from two real world problems concerning the 
Erlangian approximation of a Markovian fluid queue \cite{dl14}. 
The block-size $n$ of $\TT(U)$ is usually very large
since a bigger $n$ leads to a better Erlangian approximation, while
the size $m$ of the blocks is equal to 2 for both problems.  

We show the performances in terms of accuracy of the algorithm based
on the $\epsilon$-circulant matrix.  In Table \ref{tab:1} we report
the errors generated by Algorithm \ref{alg:ttepsc} with
$\epsilon=\cu\cdot 10^{-2}$  applied to the first problem. 
Observe that the errors are much smaller
in magnitude than $|\epsilon|$. The component-wise and norm-wise
absolute errors range around $10^{-14}-10^{-12}$, while the
componentwise relative errors deteriorate as $n$ increases; the
norm-wise relative errors moderately increase as $n$ increases.  This
behavior is expected since the use of FFT makes the algorithm stable
in norm, while the component-wise accuracy is not guaranteed.

\begin{table}
\centering
\begin{tabular}{|c|cccc|}\hline
$n$&cw-abs&cw-rel&nw-abs&nw-rel\\ \hline
128&\tt 8.4e-14  &\tt  5.8e-11 &\tt 6.5e-12  &\tt  8.0e-12 \\ \hline
256& \tt 1.1e-14   &  \tt 8.9e-11 & \tt 1.7e-12  & \tt  2.1e-12 \\ \hline
512& \tt 1.2e-14  &  \tt 2.6e-09 & \tt 7.7e-13 &  \tt  9.6e-13 \\ \hline
1024&\tt 2.2e-14   & \tt 9.1e-04&\tt 1.9e-12   &\tt  2.4e-12\\ \hline
\end{tabular}\caption{Errors generated by Algorithm \ref{alg:ttepsc}, based  on the $\epsilon$-circulant technique, with $\epsilon=\cu\cdot 10^{-2}$}\label{tab:1}
\end{table}

\begin{figure}
\begin{center}
\subfloat[Without scaling]{\includegraphics[width =6.3cm ]{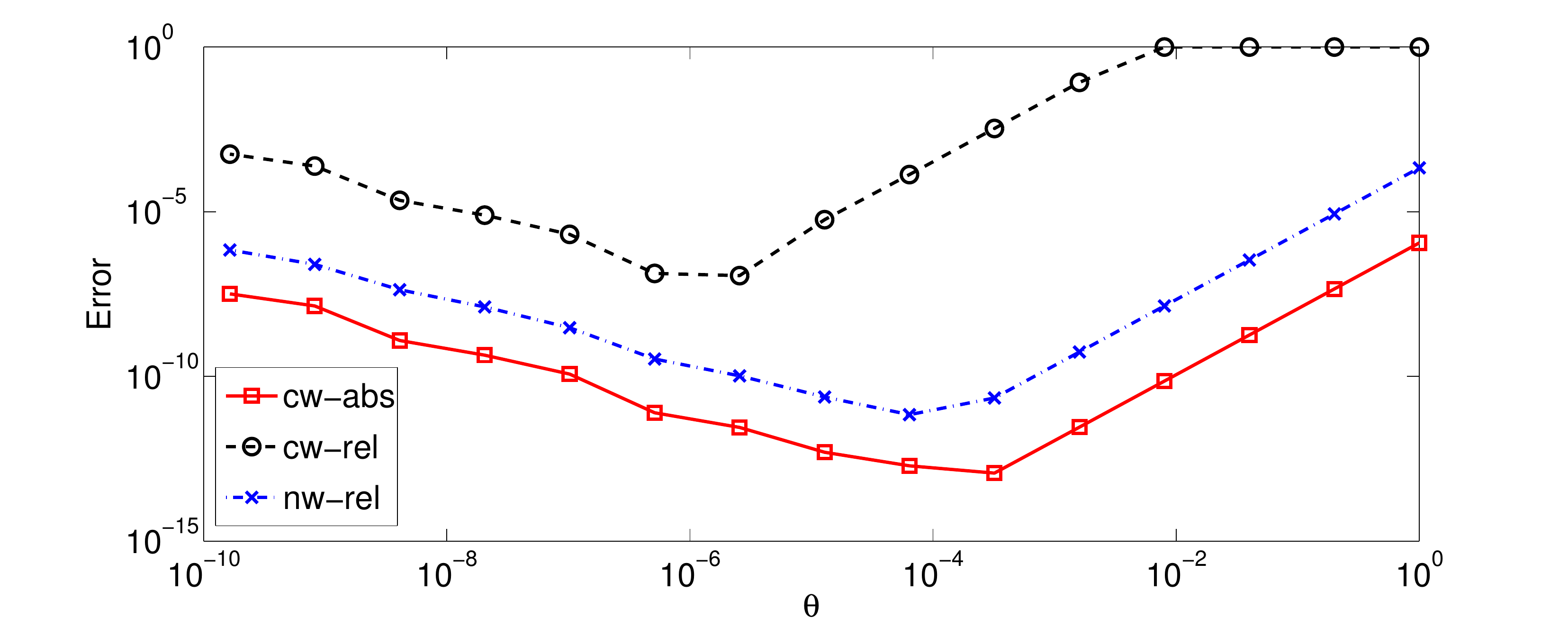}\label{fig:epca}}
\subfloat[With scaling]{\includegraphics[width = 6.3cm]{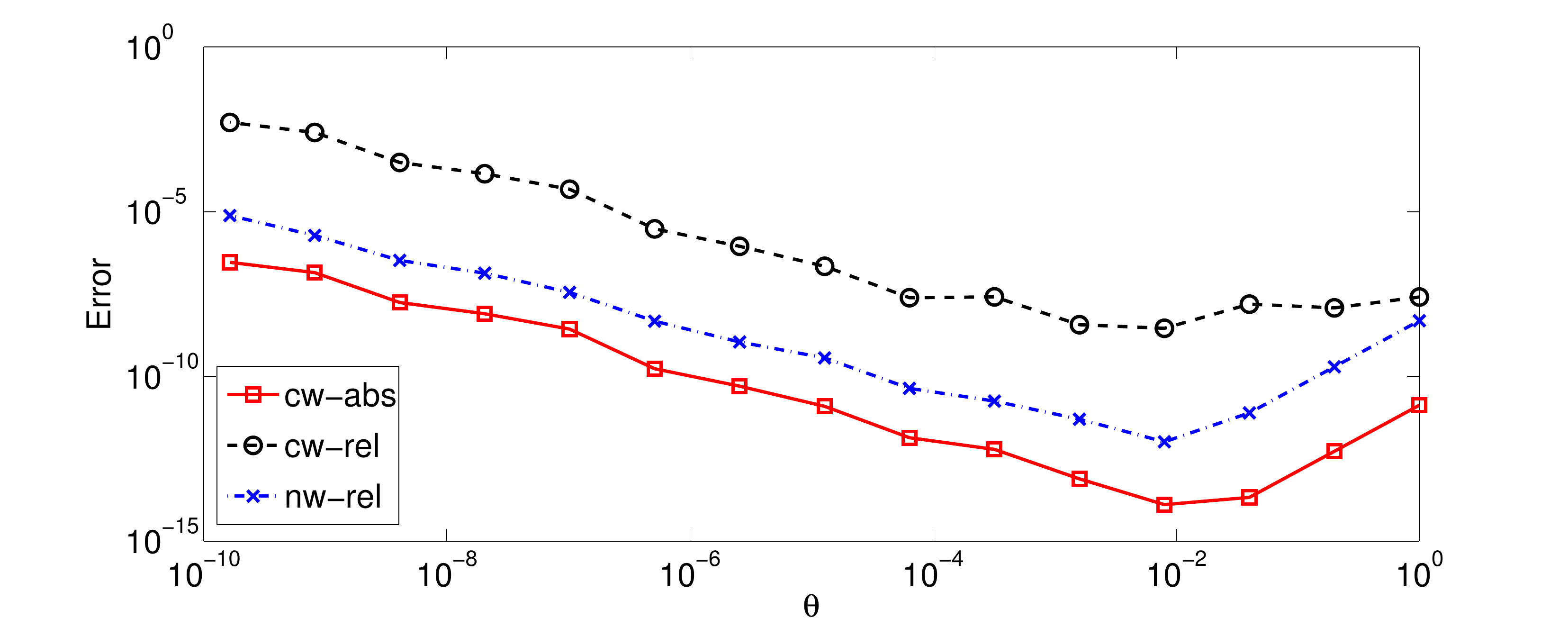}\label{fig:epcb}}
\end{center}
\caption{Error as function of $\epsilon=\cu\theta$ for the $\epsilon$-circulant algorithm, with $n=512$}\label{fig:epc}
\end{figure}

In Figure \ref{fig:epc} we report the absolute/relative component-wise
and the relative norm-wise errors as a function of
$\epsilon=\cu\cdot\theta$, with $\theta$ varying from $10^{-10}$ to
$10^0$, in the case $n=512$.  In Figure \ref{fig:epca} the scaling
technique is not applied, while in Figure \ref{fig:epcb} the scaling
is applied, as described in Algorithm \ref{alg:ttepsc}. It is worth
pointing out how the scaling allows to obtain a better accuracy, and
the best performances are obtained with a larger value of
$|\epsilon|$. Observe also that with the scaling technique the
component-wise relative error takes values close to $10^{-9}$ while
the theoretical bound is asymptotically $(2/3)\mu\approx 1.e-10$.
Another interesting remark is that the absolute component-wise errors
and the relative norm-wise errors reach a minimum value for a
moderately large value of $\theta$, and substantially increase for
values smaller than this minimum. This is due to the effect of
round-off errors, which increase as $|\epsilon|$ goes to zero.

\begin{figure}
\begin{center}
\includegraphics[width=8 cm]{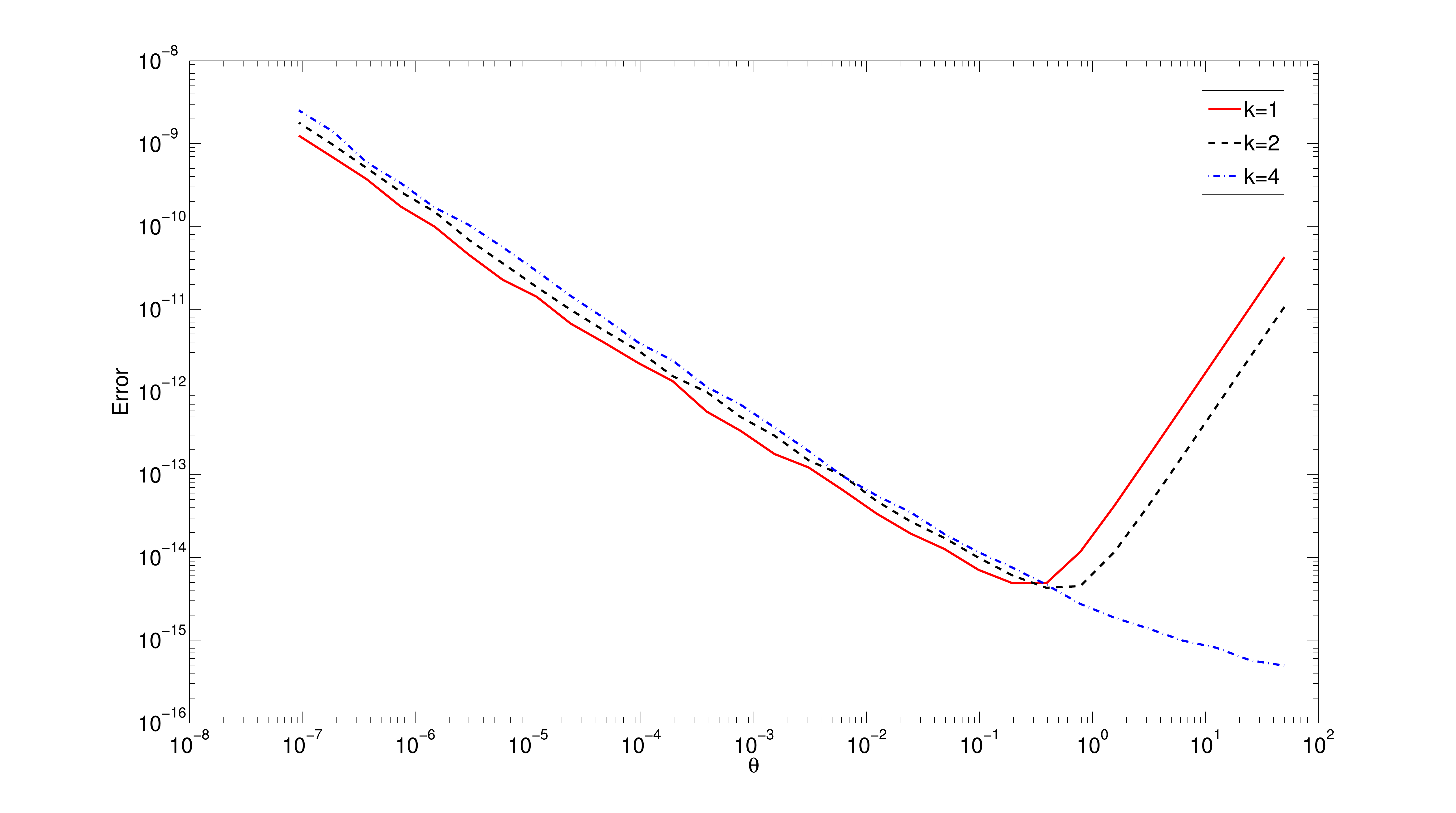}
\end{center}
\caption{Error as function of $\theta$ for the $\epsilon$-circulant algorithm, with $k$ interpolation points and $n=512$}\label{fig:epc_int}
\end{figure}

In Figure \ref{fig:epc_int} we report the normwise relative errors
obtained with the $\epsilon$-circulant technique, described in
Algorithm \ref{alg:epscirc_int}, applied with $k$ different values
$\epsilon_j=({\cu})^{1/k}\omega_k^j\theta$, for $j=0,\ldots,k-1$, where
the solution is the arithmetic mean of $e^{\mathcal
  C_{\epsilon_j}(U)}$. It is interesting to observe that using $k=2$
leads to an approximation error better than $k=1$, while for $k=4$ the
solution provided by the algorithm has an error close to the machine
precision. Actually from this picture it is possible to figure out
where the approximation errors and the roundoff errors dominate. For
$k=4$ the graph of the overall error is almost decreasing, this shows
that the approximation error is removed by the technique of averaging
the approximations obtained with different values of $\epsilon_j$.
From this behaviour one deduces that the approximation error
numerically behaves like a polynomial of degree less than 8.
This guess should be worth being investigated from a theoretical point of view.

Now consider the method based on the embedding into a circulant
matrix.  In Table \ref{tab:emb} we report the errors generated by
Algorithm \ref{alg:ttemb} with $K=4n$ applied to the first problem.
The component-wise absolute errors are of the order of the machine
precision, while the component-wise relative errors deteriorate as $n$
increases; the norm-wise relative errors remain quite small as $n$
increases.  As for the $\epsilon$-circulant method, this behavior is
expected for the use of FFT.  The accuracy of this algorithm is better
than that obtained with the $\epsilon$-circulant method.

\begin{table}
\centering
\begin{tabular}{|c|cccc|}\hline
$n$&cw-abs&cw-rel&nw-abs&nw-rel\\ \hline
128&\tt 2.0e-16  & \tt  1.4e-12&\tt 8.9e-15    &\tt 1.1e-14\\ \hline
256&\tt 4.0e-16  &\tt   2.4e-11&\tt 2.2e-14   &\tt  2.8e-14\\ \hline
512&\tt 8.5e-16 &\tt    2.5e-09&\tt 4.3e-14 &\tt   5.4e-14\\ \hline
1024&\tt 6.3e-16 &\tt 3.4e-04 &\tt 8.4e-14   &\tt  1.0e-13\\ \hline
\end{tabular}\caption{Errors generated by  Algorithm \ref{alg:ttemb}, based on the embedding technique, with $K=4n$}\label{tab:emb}
\end{table}

\begin{figure}
\centering
\subfloat[Without scaling]{\includegraphics[width =6.3cm ]{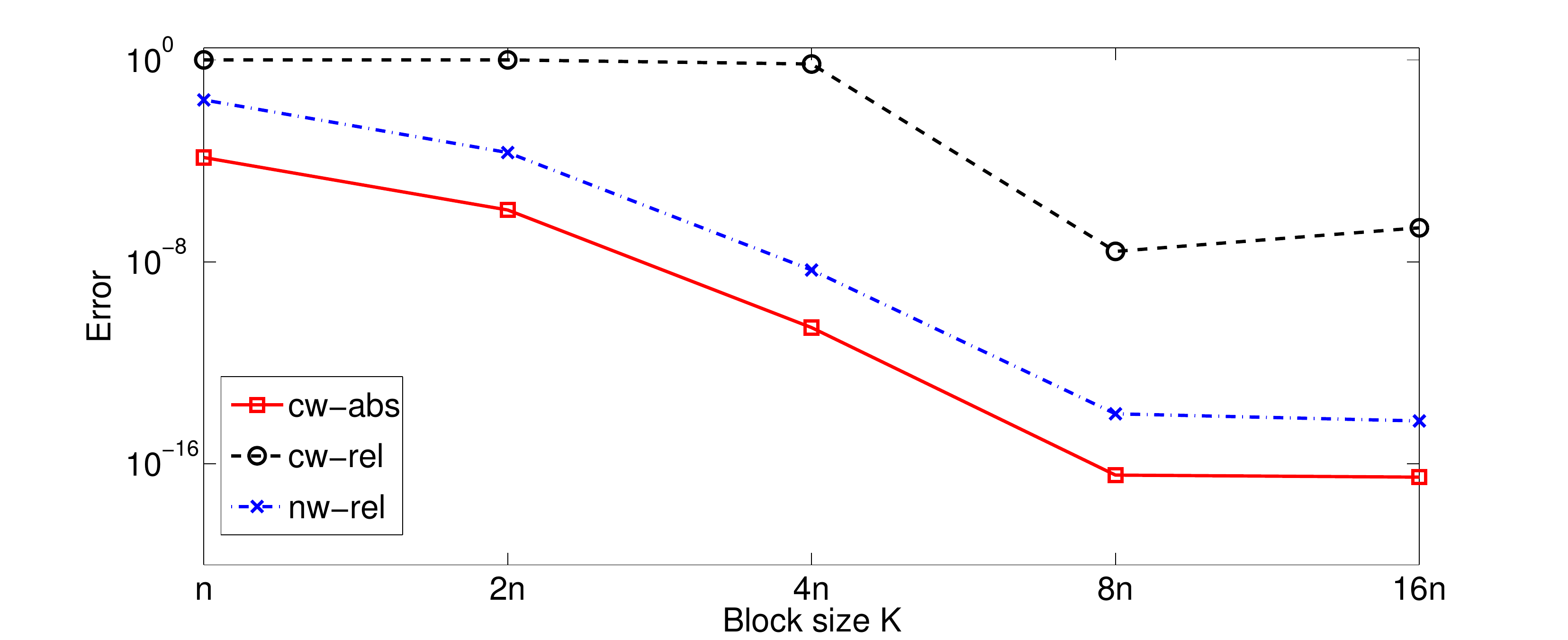}\label{fig:emba}}
\subfloat[With scaling]{\includegraphics[width = 6.3cm]{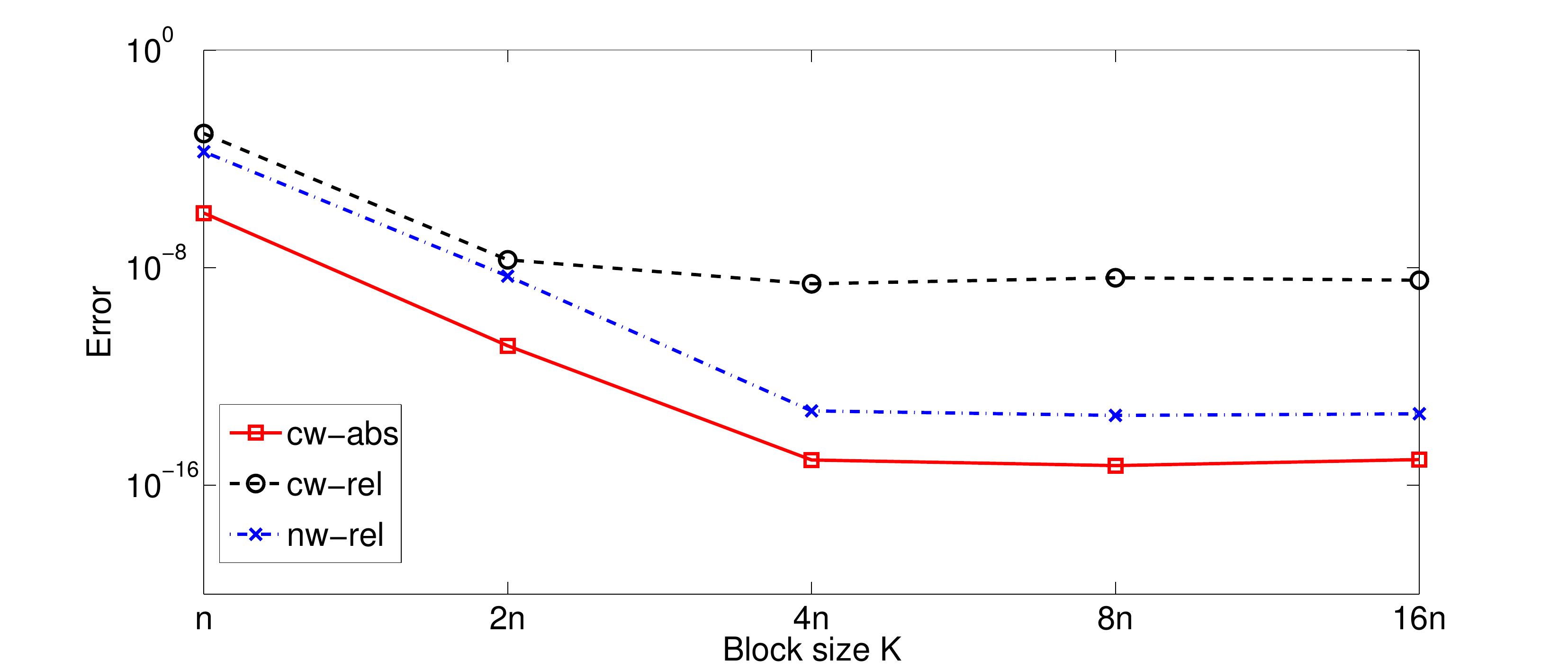}\label{fig:embb}}
\caption{Error as function of $K$ for the embedding algorithm, with $n=512$}\label{fig:emb}
\end{figure}

In Figure \ref{fig:emb} we report the absolute/relative component-wise
and relative norm-wise errors as a function of $K$, in the case
$n=512$.  In Figure \ref{fig:emba} the scaling technique is not
applied, while in Figure \ref{fig:embb} the scaling is applied, as
described in Algorithm \ref{alg:ttemb}. Also in this case it is worth
pointing out how the scaling allows to obtain a better accuracy and
optimal performances with smaller value of the block-size $K$, that is
$4n$ vs.~$8n$.

In Table \ref{tab:ttt} we report the errors generated by Algorithm
\ref{alg:ttt} based on Taylor expansion. The errors have the same
magnitude as those of Table \ref{tab:emb} for the method based on the
embedding.

\begin{table}
\centering
\begin{tabular}{|c|cccc|}\hline
$n$&cw-abs&cw-rel&nw-abs&nw-rel\\ \hline
128&\tt 4.7e-16  &\tt   9.5e-13& \tt 5.7e-15  &\tt   7.0e-15\\ \hline
256&\tt 1.8e-15   &\tt 4.3e-12 &\tt 2.2e-14   &\tt  2.8e-14\\ \hline
512& \tt 8.9e-16   &\tt  1.3e-09 &\tt 3.0e-14 &\tt  3.8e-14  \\ \hline
1024&\tt 4.8e-15  & \tt 7.7e-04&\tt 1.3e-13    &\tt 1.6e-13 \\ \hline
\end{tabular}\caption{Errors generated by  Algorithm \ref{alg:ttt}, based on Taylor expansion}\label{tab:ttt}
\end{table}

In Table \ref{tab:time} we report the CPU time in seconds, as a function of $n$,
needed by the algorithm based on $\epsilon$-circulant matrix
(\texttt{epc}), on embedding into a circulant matrix (\texttt{emb}),
on Taylor series expansion (\texttt{taylor}) and by the \texttt{expm}
function of Matlab.  The symbol ``*'' denotes an execution time greater than 100 seconds. The time needed by \texttt{expm} increases much
faster than the time needed by the other methods. The method
\texttt{epc} is the fastest, and the method based on embedding is
slightly faster than the Taylor series method when $n$ is large enough.

\begin{table}
\centering
\begin{tabular}{c|ccccc}
Algorithm~$\backslash ~ n$ & 256 & 512 & 1024 & 2048 & 4096 \\
\hline
 \texttt{epc} & 0.2 & 0.5 & 1.5 & 4.6 & 16.0\\
\texttt{emb} & 0.4 & 0.9 & 2.4 & 6.8 & 22.4 \\
\texttt{taylor} & 0.6 & 1.4 & 3.8 & 11.5 & 37.6 \\
\texttt{expm} & 0.9 & 5.9 & 327.7 & * & *
\end{tabular} 
\caption{CPU time as function of the block-size $n$}\label{tab:time}
\end{table}

\begin{figure}
\centering
\includegraphics[width =8cm ]{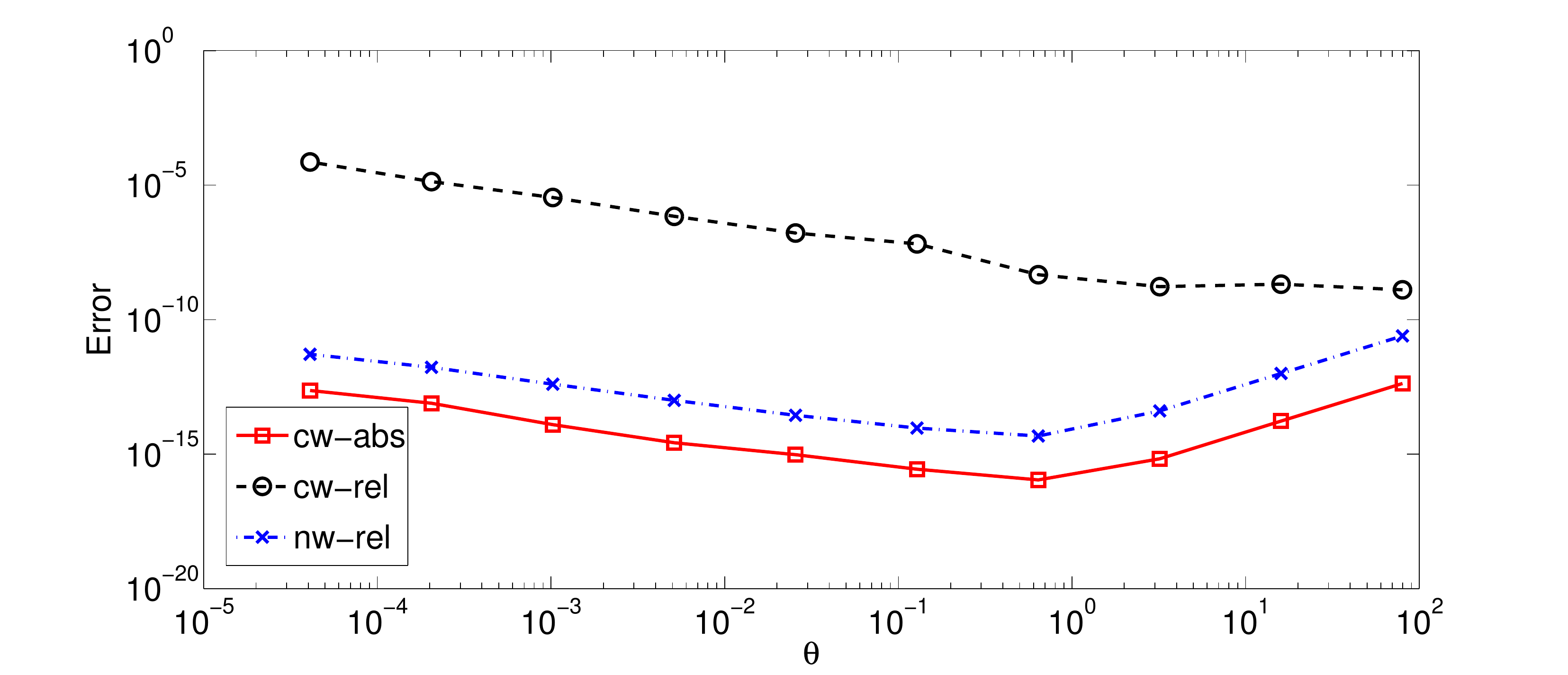}
\caption{Error as function of $\epsilon=\cu\theta$ for the $\epsilon$-circulant algorithm, with $n=512$}\label{fig:prob2.1}
\end{figure}

Concerning the second problem, we report only the results in the case
where scaling is applied. In fact, there is not much differences
between the sclaed and the unscaled versions since this problem is
already well scaled in its original formulation. In Figure
\ref{fig:prob2.1} we report the errors for the method based on
$\epsilon$-circulant matrices. It is interesting to note that the
optimal value of $|\epsilon|$ is close to 1 and that the
component-wise relative error is minimized by values of $|\epsilon|$
greater than 1. This fact, which apparently seems to be a
contradiction, is explained as follows. Large values of
$\epsilon$ generate large errors in the lower triangular part, i.e.,
the lower triangula part of $e^{\TT(U)}-e^{\CC_\epsilon(U)}$ has large
norm. On the other hand we consider the first block-row of
$e^{\CC_\epsilon(U)}$ to approximate the matrix exponential of
$\TT(U)$, therefore the errors are not influenced by a large error in
the lower triangular part.

\begin{figure}
\centering
\includegraphics[width =8cm ]{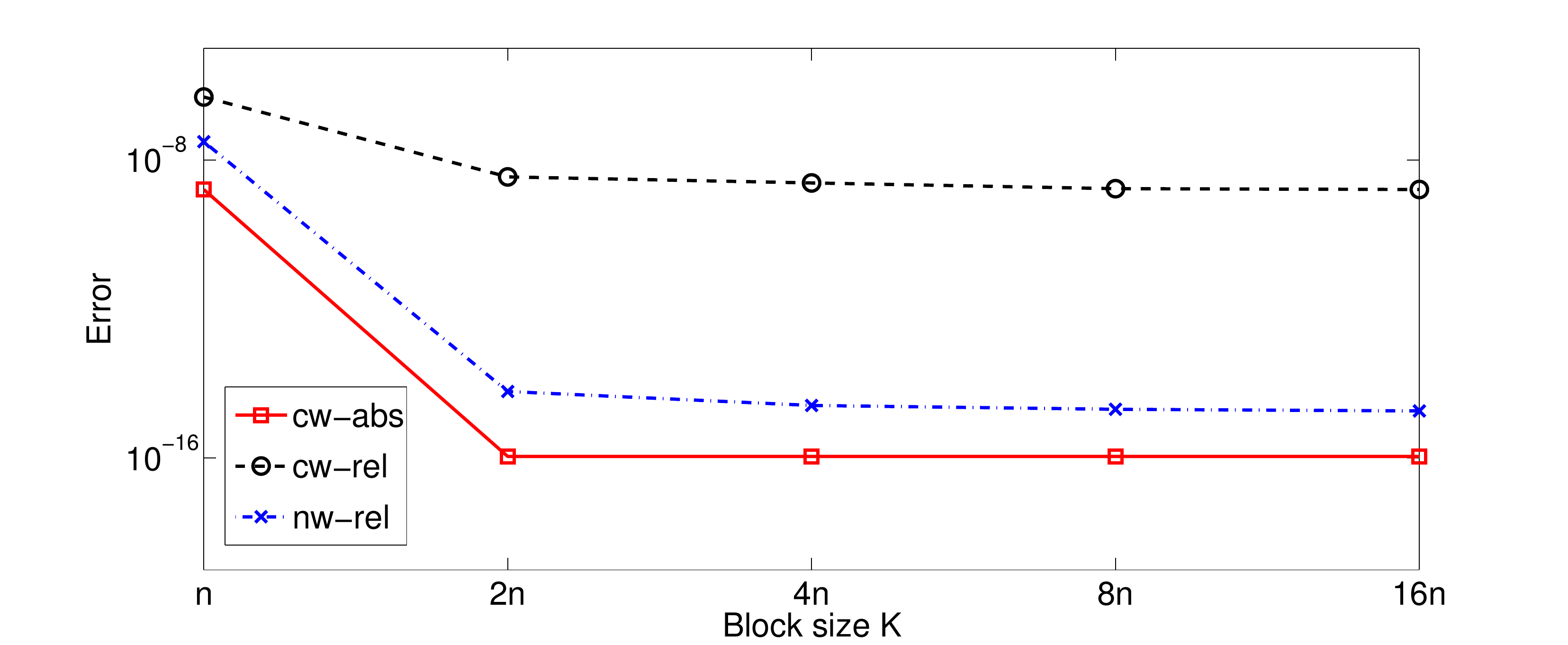}
\caption{Error as function of $K$ for the embedding algorithm, with $n=512$}\label{fig:prob2.2}
\end{figure}

In Figure \ref{fig:prob2.2} we report the errors for the algorithm based on embedding. It is relevant to observe that the errors are essentially minimized 
with an embedding of just double size.

To conclude, the method based on $\epsilon$-circulant is the fastest
one, but the accuracy of the results is lower than that provided by
the embedding and Taylor series expansion. However, by applying the averaging technique we can dramatically improve the accuracy of the $\epsilon$-circulant algorithm.

The computational time of all the structured algorithms is much lower than the cost of the general method implemented in the {\tt expm} function of Matlab
and allows to deal with matrices with huge size.

The algorithms based on embedding, on $\epsilon$-circulant matrices
are faster than the one based on Taylor series with FFT matrix
arithmetic. Moreover they are better suited for a parallel
implementation.

\end{document}